\documentclass[11pt,reqno,a4paper]{amsart_gus}
\usepackage{amssymb, amscd, amsmath, enumerate,verbatim}
\usepackage{amscd}
\usepackage[dvips,usenames]{color}

\input xy
\xyoption{all}

\numberwithin{equation}{subsection}
\newtheorem{theorem}{Theorem}[subsection]
\newtheorem{lemma}[theorem]{Lemma}
\newtheorem{proposition}[theorem]{Proposition}
\newtheorem{corollary}[theorem]{Corollary}

\theoremstyle{definition}
\newtheorem{nada}[theorem]{}
\newtheorem{definition}[theorem]{Definition}
\newtheorem{remark}[theorem]{Remark}

\newtheorem{example}[theorem]{Example}
\newtheorem{examples}[theorem]{Examples}

\newcommand{\mk}{\mathbf{k}}

\newcommand{\E}{\mathcal{E}}
\newcommand{\D}{\mathcal{D}}
\newcommand{\C}{\mathcal{C}}
\newcommand{\B}{\mathcal{B}}
\newcommand{\A}{\mathcal{A}}

\newcommand{\X}{\mathcal{X}}

\newcommand{\W}{\mathcal{W}}
\newcommand{\Scal}{\mathcal{S}}

\newcommand{\M}{\mathcal{M}}

\newcommand{\Top}{\mathbf{Top}}

\newcommand{\Cat}{\mathbf{Cat}}

\newcommand{\loc}[2]{\mathcal{#1}[\mathcal{#2}^{-1}]}

\newcommand{\Cyl}{\mathrm{Cyl}\,}
\newcommand{\Path}{\mathrm{Path}\,}
\newcommand{\Ran}{\mathrm{Ran}\,}

\newcommand{\ldf}[1]{\mathbb{L}}

\newcommand{\pChains}[1]{\mathbf{C}_{+}(#1)}

\newcommand{\MOp}{\mathbf{MOp}}
\newcommand{\Op}{\mathbf{Op}}

\newcommand{\Adgc}{\mathbf{Adgc}(\mk )}

\newcommand{\Hom}{\mathrm{Hom}\,}

\newcommand{\Diff}{{\mathbf {Diff}}}
\newcommand {\lra}{\longrightarrow}
\newcommand {\lla}{\longleftarrow}
\newcommand {\ra}{\rightarrow}

\newcommand{\PChains}[1]{\mathbf{C}_+(#1)}
\newcommand{\FPChains}[1]{\mathbf{FC}_+(#1)}

\newtheorem{Resume of notations}[theorem]{Resume of
notations}

\newtheorem{notations}[theorem]{Notation}

\newcommand{\id}{{\rm id}}
\newcommand{\Tot}{{\rm Tot}}

\begin{document}
\begin{large}

\title[A Cartan-Eilenberg approach]{A
Cartan-Eilenberg approach to homotopical algebra}

\author[F. Guill\'{e}n]{F. Guill\'{e}n }
\address[F. Guill\'{e}n  and V. Navarro]{ Departament
d'\`{A}lgebra i Geometria\\  Universitat de Barcelona\\ Gran Via 585,
08007 Barcelona (Spain)}
\author[V. Navarro] {V. Navarro }

\author[P. Pascual]{P. Pascual}
\address[P. Pascual and A. Roig]{ Departament de Matem\`{a}tica Aplicada
I\\ Universitat Polit\`{e}cnica de Catalunya\\Diagonal 647, 08028
Barcelona (Spain). }
\author[Agust\'{\i} Roig]{Agust\'{\i} Roig}
\email{fguillen@ub.edu\\vicenc.navarro@ub.edu\\pere.pascual@upc.es
\\agustin.roig@upc.es
}
\begin{abstract}
In this paper we propose an approach to homotopical algebra where
the basic ingredient is a category with two classes of
distinguished morphisms: strong and weak equivalences. These data
determine the cofibrant objects  by an extension property
analogous to the classical lifting property of projective modules.
We define a Cartan-Eilenberg category as a category with strong
and weak equivalences  such that there is  an equivalence of
categories  between its localisation  with respect to weak
equivalences and the relative localisation of the subcategory of
cofibrant objets with respect to strong equivalences. This
equivalence of categories allows us to extend the classical theory
of derived additive functors  to this non additive setting. The
main examples include Quillen model categories and  categories of
functors defined on  a category endowed with a cotriple (comonad)
and taking values on a category of complexes of an abelian
category. In the latter case there are examples in which the class
of strong equivalences is not determined by a homotopy relation.
Among other applications of our theory, we establish a very
general acyclic models theorem.

\end{abstract}

\footnotetext[1]{Partially supported by projects DGCYT MT
M2006-14575 }

\date{\today}
\maketitle
\tableofcontents

In their pioneering work \cite{CE}, H. Cartan and S. Eilenberg
defined the notion of derived functors of additive functors
between categories of modules. Their approach is based on the
characterisation of projective modules over a ring $A$ in terms of
the notions of homotopy between morphisms of complexes of
$A$-modules and quasi-isomorphisms of complexes. Projective
modules can be characterised from them: an $A$-module $P$  is
projective if  for every solid diagram
$$ \xymatrix{
& Y \ar[d]^w \\
  P \ar@{.>}[ur]^{g} \ar[r]^{f} & X  }
$$
where $w$ is a quasi-isomorphism of complexes,  and $f$ a chain
map, there is a lifting $g$ such that the resulting diagram is
homotopy commutative, and the lifting $g$ is unique up to
homotopy.

A. Grothendieck, in his Tohoku paper \cite{Gr}, introduced abelian
categories and extended Cartan-Eilenberg methods to derive
additive functors between them. Later on, Grothendieck stressed
the importance of  complexes, rather than modules, and promoted
the introduction of derived categories by J.L. Verdier.

In modern language the homotopy properties of projective complexes
can be summarised in the following manner. If $\A$ is an abelian
category with enough projective objects, then there is an
equivalence of categories
$$
\mathbf K_+(Proj(\A))\stackrel{\sim }{\lra }\mathbf
D_+(\A),\eqno{(0.1)}
$$
where $\mathbf K_+(Proj(\A))$ is the category of bounded below
chain complexes of projective objects modulo homotopy, and
$\mathbf D_+(\A)$ is the corresponding derived category. Additive
functors can therefore be derived as follows. If $F:\A\lra \B$ is
an additive functor, it induces a functor $F': \mathbf
K_+(Proj(\A))\lra \mathbf K_+(\B)$ and by the equivalence $(0.1)$,
we obtain the derived functor $\mathbb L F:\mathbf D_+(\A)\lra
\mathbf D_+(\B)$.

In order to derive non additive functors, D. Quillen, inspired by
topological methods, introduced model categories in his notes on
Homotopical Algebra \cite{Q}. Since then, Homotopical Algebra has
grown considerably as can be seen, for example, in \cite{DHKS}, \cite{Ho}, \cite{Hi}. Quillen's approach applies to
classical  homotopy theory as well as to rational homotopy,
Bousfield localisation, or more recently to simplicial sheaves or
motivic homotopy theory.

In a Quillen model category $\C$, a homotopy relation for
morphisms is defined from the axioms and one of the main results
of \cite{Q} is the equivalence
$$
\pi \C_{cf}\stackrel{\sim}{\lra }\C[\W^{-1}]\ ,\eqno{(0.2)}
$$
where $\pi \C_{cf}$ is the homotopy category of the full
subcategory $\C_{cf}$ of fibrant-cofibrant objects, and
$\C[\W^{-1}]$ is the localised category with respect to weak
equivalences. The equivalence $(0.2)$ extends the one for
projective complexes $(0.1)$ and allows derivation of functors in
this setting.

The set of axioms of model categories is, in some sense, somewhat
strong because  there are interesting categories in which to do
homotopy theory that do not satisfy all of them. Several authors
(see \cite{Br}, \cite{Ba} and others) have developed simpler
alternatives, all of them focused on laterality, asking only for a
left- (or right-) handed version of Quillen's set of axioms. All
these alternatives are very close to Quillen's formulation.

Here we propose another approach which is closer to the original
development by Cartan-Eilenberg. The initial data are two classes of
morphisms $\mathcal S$ and $\W$ in a category $\C$, with $\mathcal
S\subset \overline{\W}$, which we call strong and weak equivalences,
respectively. We define an object $M$ of $\C$ to be cofibrant if for
every solid diagram
$$ \xymatrix{
& Y \ar[d]^w \\
 M \ar@{.>}[ur]^{g} \ar[r]^{f} & X\ ,  }
$$
where $w$ is a weak equivalence and $f:M\lra X$ is a morphism in
$\C$, there is a unique lifting $g$ in $\C[{\mathcal S}^{-1}]$
such that the diagram is commutative in $\C[\mathcal S^{-1}]$.
 We say that $\C$ is a Cartan-Eilenberg
category if  it has enough cofibrant objects, that is, if  each
object $X$ in $\C$ is isomorphic in $\C[\W^{-1}] $ to a cofibrant
object.  In that case the functor
$$
\C_{cof}[\mathcal S^{-1},\C]\stackrel{\sim}{\lra }\C[\mathcal
W^{-1}]\eqno{(0.3)}
$$
is an equivalence of categories, where $\C_{cof}[\mathcal
S^{-1},\C]$ is the full subcategory of $\C[\mathcal S^{-1}]$ whose
objects are the cofibrant objects of $\C$.

In a Cartan-Eilenberg category we can derive functors exactly in
the same way as Cartan Eilenberg. If $\C$ is a Cartan-Eilenberg
category and $F:\C\lra \D$ is a functor which sends strong
equivalences to isomorphisms, $F$ induces a functor $ F':
\C_{cof}[\mathcal S^{-1},\C]\lra \D$ and by the equivalence
$(0.3)$, we obtain the derived functor $\mathbb L F:
\C[\W^{-1}]\lra \D$.

Each  Quillen model category produces a Cartan-Eilenberg category:
the category of its fibrant objects, with $\mathcal S$ the class
of left homotopy equivalences and $\W$ the class of weak
equivalences. Nevertheless, note the following differences with
Quillen's theory. First, in the Quillen context the class
$\mathcal S$ appears as a consequence of the axioms while
fibrant/cofibrant objects are part of them. Second, cofibrant
objects in our setting are homotopy invariant, in contrast with
cofibrant objects in Quillen model categories. Actually, in a
Quillen category of fibrant objects, an object is Cartan-Eilenberg
cofibrant if and only if it is homotopy equivalent to a Quillen
cofibrant one.

Another example covered by our presentation is that of Sullivan's
minimal models. We define minimal objects in a Cartan-Eilenberg
category, and call it a Sullivan category, if any object has a
minimal model. As an example, we interpret some results of
\cite{GNPR1} as saying that the category of modular operads over a
field of characteristic zero is a Sullivan category.

In closing this introduction,  we want to highlight  the
definition of Cartan-Eilenberg structures coming from a cotriple.
If $\X$ is a category with a cotriple $G$,  $\A$ is an abelian
category and $\mathbf{C}_+(\A)$ denotes the category of bounded
below chain complexes of $\A$, we define a structure of
Cartan-Eilenberg category on the functor category
$\mathbf{Cat}(\X, \mathbf{C}_+(\A))$ (see Theorem
\ref{modelsaciclics1}). We apply this result to obtain theorems of
the acyclic models kind, extending results in \cite{B} and
\cite{GNPR2}. We stress that in these examples the class of strong
equivalences $\mathcal S$ does not come from a homotopy relation.

\emph{Acknowledgements}. We thank C. Casacuberta, B. Kahn and G.
Maltsiniotis for their comments on an early draft of this paper.
We are also indebted to the referee for his kind remarks and
critical observations.

\section{Localisation of Categories}

In this section we collect for further reference some mostly
well-known facts about localisation of  categories, and we
introduce the notion of relative localisation of a subcategory,
which plays an important role in the sequel.

\subsection{Categories with weak equivalences}

\begin{nada}
By a {\em category with weak equivalences} we understand a pair
$(\C,\W)$ where $\C$ is a category and $\W$ is a class of
morphisms of $\C$. Morphisms in $\W$ will be called \emph{weak
equivalences}.

We always assume that $\W$ is stable by composition and contains
all the isomorphisms of $\C$, so that we can identify $\W$ with a
subcategory of $\C$.
\end{nada}

\begin{nada}
Recall that the \emph{category of fractions}, or
\emph{localisation, of $\C$ with respect to }$\W$ is a category
$\loc{C}{W}$ together with a functor $\gamma:\C\longrightarrow
\mathcal C\left[\mathcal W^{-1}\right]$ such that:

\begin{itemize}
\item[(i)] For all $w\in \W$, $\gamma(w)$ is an isomorphism.

\item[(ii)] For any category $\D$ and any functor $F:\C\lra \D$
that transforms morphisms $w\in \W$ into isomorphisms, there
exists a unique functor $F':\C[\W^{-1}]\lra \D$ such that $F'\circ
\gamma =F$.
\end{itemize}
 The uniqueness condition on $F'$ implies immediately that, when it
exists, the localisation is uniquely defined up to isomorphism.
The localisation  exists if $\W$ is small and, in general, the
localisation always exists in a higher universe.
\end{nada}

\begin{nada} We say that the class of weak equivalences $\W$ is
\emph{saturated} if a morphism $f$ of $\C$ is in $\W$ when $\gamma
f$ is an isomorphism. The \emph{saturation} $\overline{\W}$ of
$\W$ is the pre-image by $\gamma$ of the isomorphisms of
$\C[\W^{-1}]$. It is the smallest saturated class of morphisms of
$\C$ which contains $\W$. Maybe it is worth pointing out that we
do \emph{not} assume that $\W$ verifies the usual \emph{2 out of 3
property}. In any case, the saturation $\overline{\W}$ always
does.
\end{nada}

\subsection{Hammocks}
 We describe the localisation of
categories by using Dwyer-Kan hammocks (\cite{DK}).  Given a
category with weak equivalences $(\C,\W)$ and two objects X and Y
in $\C$, a {\em $\W$-zigzag} $f$ from $X$ to $Y$ is a finite
sequence of morphisms of $\C$, going in either direction, between
$X$ and $Y$,
$$
\xymatrix{ f: \; X  \ar@{-}[r] & \bullet\ar@{-}[r] &\bullet \dots
\bullet \ar@{-}[r] &\bullet\ar@{-}[r]  & Y   },
$$
where the morphisms going from right to left are in $\W$. We call the number of morphisms in the sequence the \emph{length} of the $\W$-zigzag. Because each $\W$-zigzag is a diagram, it has a
\emph{type}, its index category. A morphism from a $\W$-zigzag $f$
to a $\W$-zigzag $g$ of the same type is a commutative diagram in
$\C$,
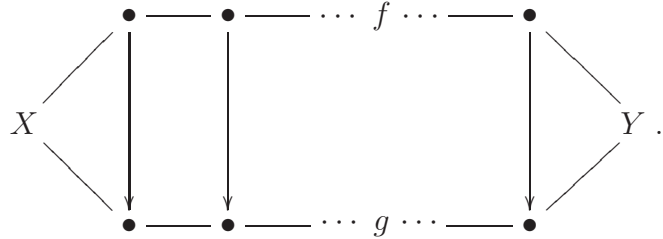
$$
\xymatrix{
    & {\bullet} \ar@{-}[r] \ar[dd]& {\bullet} \ar[dd] \ar@{-}[r] &
    \cdots \ f\ {\cdots} & {\bullet} \ar[dd]\ar@{-}[l]
    \ar@{-}[dr]&             \\
 X  \ar@{-}[ur] \ar@{-}[dr]    &  &  &  &  &       Y \; .     \\
    & {\bullet}   \ar@{-}[r]  & {\bullet}   \ar@{-}[r]  & \cdots\ g\
    {\cdots} &
    {\bullet}  \ar@{-}[l] \ar@{-}[ur] &
}
$$

A hammock  between two $\W$-zigzags $f$ and $g$ from $X$ to $Y$ of
the same type is a finite sequence of morphisms of zigzags going
in either direction. More precisely, it is a commutative diagram
$H$ in $\C$
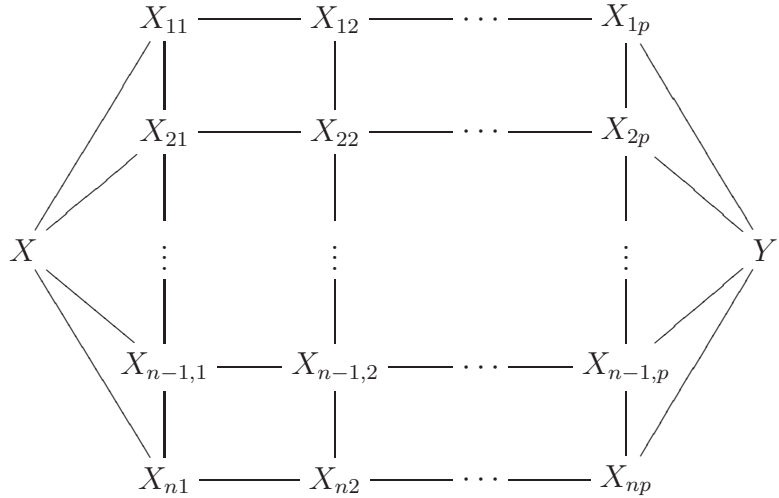
$$
\xymatrix{
    & {X_{11}} & {X_{12}} \ar@{-}[l] \ar@{-}[r] & {\cdots} & {X_{1p}}
    \ar@{-}[l]
    \ar@{-}[ddr]&             \\
    & {X_{21}} \ar@{-}[u] \ar@{-}[d] & {X_{22}} \ar@{-}[l] \ar@{-}[u]
     \ar@{-}[r] \ar@{-}[d] & {\cdots} & {X_{2p}}
    \ar@{-}[dr] \ar@{-}[u] \ar@{-}[l] \ar@{-}[d] &             \\
X \ar@{-}[uur] \ar@{-}[ur] \ar@{-}[dr] \ar@{-}[ddr]
    & {\vdots} & {\vdots} &     & {\vdots} &      Y      \\
    & {X_{n-1,1}} \ar@{-}[d]  \ar@{-}[u] & {X_{n-1,2}}
    \ar@{-}[l] \ar@{-}[d] \ar@{-}[r] \ar@{-}[u] & {\cdots} & {X_{n-1,p}}
    \ar@{-}[ur] \ar@{-}[d] \ar@{-}[l] \ar@{-}[u] &             \\
    & {X_{n1}} & {X_{n2}} \ar@{-}[l] \ar@{-}[r] & {\cdots} & {X_{np}}
    \ar@{-}[uur] \ar@{-}[l] &             \\
}
$$
such that

\begin{enumerate}[(i)]
\item  in each column of arrows, all (horizontal) maps go in the
same direction, and if they go to the left they are in $\W$ (in
particular, any row is a $\W$-zigzag),

\item in each row of arrows, all (vertical) maps go in the same
direction, and they are arbitrary maps in $\C$,

\item  the top $\W$-zigzag is $f$ and the bottom is $g$.

\end{enumerate}

If there is a hammock $H$ between $f$ and $g$, and $f'$ is a
$\W$-zigzag obtained from $f$ adding identities, then adding the
same identities in the  hammock  $H$ and in the $\W$-zigzag $g$ we
obtain a new $\W$-zigzag $g'$ and a    hammock  $H'$  between $f'$
and $g'$.

We say that two $\W$-zigzags $f,g$ between $X$ and $Y$ are
\emph{related} if there exist $\W$-zigzags $f'$ and $g'$ of the
same type, obtained from $f$ and $g$ by adding identities, and a
hammock $H$ between $f'$ and $g'$. This is an equivalence relation
between $\W$-zigzags. For instance, if in a $\W$-zigzag $f$ there
exist two consecutive arrows in the same direction, then $f$ is
equivalent to the $\W$-zigzag obtained from $f$ composing these
two arrows, as follows from the following diagram
$$
\xymatrix{ X_1\ar[d]^{\id}\ar[r]^ {f_1}&X_2\ar[d]^{f_2}\ar[r]^{f_2}
&X_3\ar[d]^{\id}\\
X_1\ar[r]^{f_2f_1}&X_3\ar[r]^{\id}&X_3 \; .}
$$

Furthermore, since $\W$ is closed by composition and contains the
isomorphisms, we can add identities, if necessary, and compose two
consecutive arrows in the same direction in such a way that each
$\W$-zigzag $f$ is related to a $\W$-zigzag of the form
$$
\xymatrix{  X  \ar[r] & \bullet &\bullet\ar[l] \ar[r]&\bullet\dots
\bullet &\ar[l]\bullet\ar[r]  & Y   },
$$
that is, two consecutive morphisms always go in opposite
directions and the first and the last morphisms  go to the right.
One such $\W$-zigzag will be called an \emph{alternating
$\W$-zigzag}.

Let $\C_\W$ be the category whose objects are the objects of $\C$
where, for any two objects $X$, $Y$, the morphisms from $X$ to $Y$
are the equivalence classes of $\W$-zigzags from $X$ to $Y$, with
composition being the juxtaposition of $\W$-zigzags.

\begin{theorem} \emph{(\cite{DHKS}, 33.10).}  The category  $\C_\W$,
together with the obvious functor $\C\lra \C_\W$ is a solution to
the universal problem of the category of fractions
$\C\left[\W^{-1}\right]$.
\end{theorem}
In the cited reference there is a general hypothesis which
concerns the class $\W$, which is not necessary for this result.

\begin{nada}
The localisation functor $\gamma:\C\lra \C[\W^{-1}]$ induces a
bijective map on the class of objects. In order to simplify the
notation, if $X$ is an object of $\C$, sometimes we will use the
same letter $X$ to denote its  image $\gamma(X)$ in the localised
category
  $\C[\W^{-1}]$.

We denote by $\Cat_\W(\C,\D)$ the category of functors from $\C$
to $\D$ that send morphisms in $\W$ to isomorphisms. The
definition of the category of fractions means that for any
category $\D$, the functor
$$
\gamma^*:\Cat(\C[\W^{-1}],\D)\lra \Cat_\W(\C,\D),\quad G\mapsto
G\circ \gamma
$$
induces a bijection on the class of  objects. From the previous
description of the localised category we deduce that  $\gamma^*$
is an isomorphism of categories. In particular, the functor
$$\gamma^*:\Cat(\C[\W^{-1}],\D)\lra \Cat(\C,\D)
$$
is fully faithful.

\end{nada}

\subsection{Categories with a congruence} There are some situations
where it is possible to give an easier presentation of morphisms
of the category $\C[\W^{-1}]$, for example, when there is a
calculus of fractions (see \cite{GZ}). In this section we present
an even simpler situation which will occur later, namely the localisation
provided by some quotient categories.

\begin{nada}
Let $\C$ be a category and \; $\sim$ \; a \emph{ congruence} on
$\C$, that is, an equivalence relation between morphisms of $\C$
which is compatible with composition (\cite{ML}, page 51). We
denote by $\C/\!\!\sim$ the  \emph{quotient category}, and by $\pi
:\C \longrightarrow \C /\!\! \sim$ the universal canonical
functor. We denote by $\mathcal S$ the class of  morphisms
$f:X\lra Y$ for which there exists a morphism $g: Y
\longrightarrow X$ such that $fg \sim 1_Y$ and $gf \sim 1_X$. We
will call $\mathcal S$ the class of equivalences associated to
$\sim$.
\end{nada}

\begin{nada}
If  $\sim $ is a congruence, in addition to the quotient category
$\C/\sim$, one can also consider the localised category
$\delta:\C\lra \C[\mathcal S^{-1}]$ of $\C$ with respect to the
class $\mathcal S$ of equivalences defined by this congruence. We
study when they are equivalent.
\end{nada}

\begin{proposition}\label{congruencia}
Let $\sim$ be a congruence and $\mathcal{S}$ the associated class
of equivalences. If $\mathcal{S}$ and $\sim$ are compatible, that
is, if $f \sim g$ implies $\delta f = \delta g$, then the
categories $\C/\!\!\sim$ and $\C[\mathcal{S}^{-1}]$ are
canonically isomorphic.
\end{proposition}
\begin{proof}
If $\mathcal{S}$ and $\sim$ are compatible, the canonical functor
$\delta:\C\lra \C[\mathcal S^{-1}]$  induces a functor
$\phi:\C/\sim \lra \C[\mathcal S^{-1}]$ such that
$\phi\circ\pi=\delta$. Therefore, any functor $F:\C\lra \D$ which
sends morphisms in $\mathcal S$ to isomorphisms factors in a
unique way through $\pi$, hence $\pi:\C\lra \C/\sim$ has the
universal property of localisation.
\end{proof}

\begin{example}\label{cilindro}
The congruence $\sim$ is compatible with its class $\mathcal S$ of
equivalences when it may be expressed by a cylinder object, or
dually by a path object.

Given $X\in \rm{Ob}\ \C$,  a \emph{cylinder object over} $X$ is an
object $\Cyl(X)$ in $\C$ together with morphisms $i_0,i_1:X\lra
\Cyl(X)$ and $p:\Cyl(X)\lra X$ such that $p\in \mathcal S$ and
$p\circ i_0=\id_X=p\circ i_1$.

Now, suppose that the congruence is determined by   cylinder
objects in the following way:

``Given $f_0,f_1:X\lra Y$, $f_0 \sim f_1 $ if and only if there
exists a morphism $H: \Cyl(X) \rightarrow Y$ such that $Hi_0 =
f_0$ and $Hi_1 = f_1$".

 Then
 $\sim$ and $\mathcal{S}$ are
compatible. In fact, if $f_0 \sim f_1 $, then  we have the
$\mathcal S$-hammock

$$
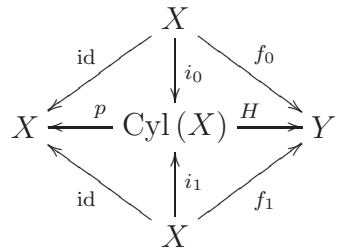
\xymatrix{
                 & X \ar[dl]_{\id} \ar[d]^{i_0} \ar[dr]^{f_0} &      \\
X    & {\Cyl (X)} \ar[l]_p \ar[r]^H           &  Y   \\
                 & X\ar[ul]^{\id} \ar[u]_{i_1} \ar[ur]_{f_1}
                      &       }
$$

between $f_0$ and $f_1$, which shows that $\delta (f_0) = \delta
(f_1)$ in $\C[\mathcal{S}^{-1}]$.

More generally,
 $\sim$ and $\mathcal{S}$ are
compatible if $\sim$ is the equivalence relation transitively
generated by a cylinder object.
\end{example}

\subsection{Relative localisation of a subcategory}
Let $\sim$ be a congruence on a category $\C$. If $i:\M\lra \C$ is
a full subcategory, there is an induced congruence on $\M$ and the
quotient category $\M/\sim$ is a full subcategory of $\C/\sim$.
Nevertheless, if $\mathcal S$ denotes the class of equivalences
associated to $\sim$, and $\mathcal S_\M$ the morphisms in $\M$
which are in $\mathcal S$, the functor $\overline{i}:\M[\mathcal
S_\M^{-1}]\lra \C[\mathcal S^{-1}]$ is not faithful, in general.
More generally, if $\E$ is an arbitrary class of morphisms in
$\C$, the functor $\overline{i}: \M[\E_\M^{-1}]\lra \C[\E^{-1}]$
is neither faithful nor full.

To simplify the notation, in the situation above we write
$\M[\E^{-1}]$ for $\M[\E_\M^{-1}]$.

\begin{definition}\label{localizacionrelativa}
Let $(\C, \E)$ be a category with weak equivalences and $\M$ a
full subcategory. The \emph{relative localisation of the
subcategory}  $\M$ of $\C$ with respect to $\E$, denoted by
$\M[\E^{-1},\C]$, is the full subcategory of $\C[\E^{-1}]$ whose
objects are those of $\M$.
\end{definition}

This relative localisation is necessary in order to express the
main results of this paper (e.g. Theorem \ref{CEimplies models}).
In  Remark \ref{minimal no localiza} we will see an interesting
example where the relative localisation $\M[\E^{-1},\C]$ is not
equivalent to the localisation $\M[\E^{-1}]$. However, in some
common situations there is no distinction between them, as for
example in the proposition below, which is an abstract generalised
version of Theorem III.2.10 in \cite{GMa}.

\begin{proposition}
Let $(\C, \E)$ be a category with weak equivalences and $\M$ a
full subcategory. Suppose that $\E$ has a right calculus of
fractions and that for every morphism $w:X\lra M$ in
$\overline{\E}$, with $M\in Ob\, \M$, there exists a morphism
$N\lra X$ in $\overline{\E}$, where $N\in Ob\, \M$. Then
$\overline{i}:\M[\E^{-1}]\lra \M[\E^{-1},\C]$ is an equivalence of
categories.
\end{proposition}
\begin{proof}
Let's prove that $\overline{i}$ is full: if $f=g
\sigma^{-1}:M_1\lla X\lra M_2$ is a morphism in $\C[\E^{-1}]$
between objects of $\M$, where $\sigma\in \E$, take a weak
equivalence $\rho:N\lra X$ with $N\in Ob\ \M$, whose existence is
guaranteed by hypothesis. Then $f=g \rho (\sigma\rho)^{-1}$ is a
morphism of $\M[\E^{-1}]$. The faithfulness is proved in a similar
way.
\end{proof}

\section{Cartan-Eilenberg categories}

In this section we define cofibrant objects in a relative setting
given by two classes of morphisms, as a generalisation of
projective complexes in an abelian category. Then we  introduce
Cartan-Eilenberg categories and give some criteria to prove that a
given category is Cartan-Eilenberg. We also relate these notions
with Adams' study of localisation in homotopy theory, \cite{A}.

\subsection{Models in a category with strong and weak equivalences}

Let $\C$ be a category and $\mathcal {S}, \W$ two classes of
morphisms of $\C$. Recall that our classes of morphisms are closed
under composition and contain all isomorphisms, but, generally
speaking, they are not saturated.

\begin{definition}\label{def(C,S,W)}
We say that $(\C,\mathcal S,\W)$ is a {\it category with strong
and weak equivalences} if $\mathcal{S}\subset\overline{\W}$.
Morphisms in $\mathcal S$ are called \emph{strong equivalences}
and those in $\W$ are called \emph{weak equivalences}.
\end{definition}

The basic example of category with strong and weak equivalences is
the category of bounded below chain complexes of $A$-modules
$\PChains{A}$, for a commutative ring $A$, with $\mathcal S$ the
class of homotopy equivalences and $\W$ the class of
quasi-isomorphisms.

\begin{notations}\label{notations2}
It is convenient to fix some notation for the rest of the paper.
 Let $(\C,\mathcal S, \W)$ be a
category with strong and weak equivalences. We denote by
$\delta:\C\lra \loc{C}{S}$ and $\gamma:\C\lra \C[\W^{-1}]$ the
canonical functors. Since $\mathcal S\subset \overline{\W}$, the
functor $\gamma $ factors through $\delta$ in the form
$$
\xymatrix{ {\C} \ar[rrr]^-{\gamma} \ar[rrd]_{\delta}& &
&  {\loc{C}{W} \cong \loc{C}{S}[\delta(\W)^{-1}]}\, .  \\
     && \loc{C}{S} \ar[ru]_{\gamma'}   &
}
$$
\end{notations}

\begin{definition}\label{defmodel2}
Let $(\C,\mathcal S, \W)$ be a category with strong and weak
equivalences, $\M$ a full subcategory of $\C$ and $X$ an object of
$\C$.  A \emph{left $(\mathcal S,\W)$-model} of $X$, or simply a
\emph {left model},  in $\M$ is an object $M$ in $\M$ together
with a morphism $\varepsilon:M\lra X$ in $\C[\mathcal S^{-1}]$
which is an isomorphism in $\C[\W^{-1}]$.

We say that there are \emph{enough left models} in $\M$, or that
$\M$ is a \emph{subcategory of left models} of $\C$, if each
object of $\C$ has a left model in $\M$.
\end{definition}

\subsection{Cofibrant objects}

\begin{definition}\label{defcofobj}
Let  $(\C,\mathcal S,\W) $ be a category with strong and weak
equivalences. An object $M$ of $\C$ is called \emph{$(\mathcal
S,\W)$-cofibrant}, or simply \emph{cofibrant}, if for each
morphism $w:Y\lra X$ of $\C$ which is in $\W$ the map
$$
w_*:\loc{C}{S}(M, Y)\longrightarrow \loc{C}{S}(M, X),\; g\mapsto
w\circ g
$$
is bijective.
\end{definition}

That is to say, cofibrant objects are defined by a lifting
property, in $\loc{C}{S}$, with respect to weak equivalences: for
any solid-arrow diagram such as
$$ \xymatrix{
& Y \ar[d]^w \\
  M \ar@{.>}[ur]^{g} \ar[r]^{f} & X  }
$$
with $w\in\W$ and $f\in\loc{C}{S}(M,X)$, there exists a
\emph{unique} morphism $g\in\loc{C}{S}(M,Y)$ making  the triangle
commutative in $\loc{C}{S}$.

\begin{proposition}\label{retracto}
 Every retract of a cofibrant object is cofibrant.
\end{proposition}
\begin{proof}
If $N$ is a retract of a cofibrant object $M$ and $w:Y\lra X$ is a
weak equivalence, the map $w^N_*:\loc{C}{S}(N,Y)\lra
\loc{C}{S}(N,X)$ is a retract of the bijective map
$w^M_*:\loc{C}{S}(M,Y)\lra \loc{C}{S}(M,X)$, hence it is also
bijective. Therefore $N$ is cofibrant.
\end{proof}

Cofibrant objects are characterised as follows (cf. \cite{Sp},
Proposition 1.4).

\begin{theorem}\label{lemallave}
Let  $(\C,\mathcal S,\W) $ be a category with strong and weak
equivalences, and $M$ an object of $\C$. The following conditions
are  equivalent.
\begin{itemize}
\item[(i)] $M$ is cofibrant.

\item[(ii)]  For each $X\in Ob\; \C$, the map
$\gamma'_X:\loc{C}{S}(M,X)\longrightarrow {\C[\W^{-1}]} ( M, X)$
is bijective.
\end{itemize}
\end{theorem}

\begin{proof} Firstly, let us see that (i) implies  (ii).
First of all, if $M$ is cofibrant,  the functor $$F:\C[\mathcal
S^{-1}]\lra \mathbf {Sets}, \; X\mapsto \C[\mathcal S^{-1}](M,X)$$
sends morphisms in $\delta(\W)$ to isomorphisms in $\mathbf
{Sets}$. Therefore this functor induces a functor on the
localisation
$$
F':\C[\W^{-1}]\lra \mathbf{Sets}
$$
such that $F'(\gamma'(f))=F(f)$ for each  $f\in\C[\mathcal
S^{-1}](X,Y)$. In addition, $\gamma'$ induces a natural
transformation
$$
\gamma':F'\lra \C[\W^{-1}](M,-)\,.
$$

Let $X$ be an object of $\C$. To see that
$$\gamma'_X:F'(X)=\loc{C}{S}(M,X)\longrightarrow {\C[\W^{-1}]} (
M, X)$$ is bijective we define a map
$$\Phi:{\C[\W^{-1}]} ( M, X)\lra F'(X)$$
which is inverse of $\gamma'_X$. Let $f\in \C[\W^{-1}](M,X)$,
then, since $F'$ is a functor, we have a map $$F'(f):F'(M)\lra
F'(X).$$ We define $\Phi(f):=F'(f)(\id_{M})$.

By the commutativity of the diagram
$$
\xymatrix{ F'(M)\ar[dd]^{\gamma_M'}\ar[rr]^{F'(f)}&&F'(X)
\ar[dd]^{\gamma'_X}\\
\\ \C[\W^{-1}](M,M)\ar[rr]^{f_*}&&\C[\W^{-1}](M,X) }
$$
we obtain
$$\gamma'_X(\Phi(f))=\gamma'_X(F'(f)(\id_M))=f_*(\gamma'_M(\id_M))=f.$$
Also, given a morphism $g\in \C[\mathcal S^{-1}](M, X)$, we have
$$\Phi(\gamma'_X(g))=F'(\gamma'_X(g))(\id_M)=F(g)(\id_M)=g,$$ so
$\Phi$ is the inverse of $\gamma '_X$, thus we obtain (ii).

Next,  (i) follows   from  (ii), since, if (ii) is satisfied, for
each $w\in \C(Y, X)$ which is in $ \W$, we have a commutative
diagram
$$\xymatrix{
 \loc{C}{S}(M,Y)\ar[dd]^{w_*}\ar[rr]^{\gamma'_Y}_\cong
 &&\C[\W^{-1}]( M, Y)\ar[dd]^{w_*}_{\cong}\\
&&\\
 \loc{C}{S}(M,X)\ar[rr]^{\gamma'_X}_{\cong}&&\C[\W^{-1}]( M, X)
 }
$$
where  three  of the arrows are bijective; thus, so is the fourth.
\end{proof}

\begin{nada}\label{+notacions}
We denote by $\C_{cof}$ the full  subcategory of $\C$ whose
objects are the cofibrant objects of $\C$,  by
$$i:\C_{cof}[\mathcal{S}^{-1},\C]\lra \loc{C}{S}$$ the inclusion
functor, and by $$j:\C_{cof}[\mathcal{S}^{-1},\C]\lra
\C[\W^{-1}]$$ the composition $j:=\gamma'\circ i$.

>From Definition \ref{defcofobj}, it follows that an object
isomorphic in $\C[\mathcal S^{-1}]$ to a cofibrant object  is also
a cofibrant object, therefore $\C_{cof}[\mathcal{S}^{-1},\C]$ is a
replete subcategory of $\loc{C}{S}$. (We recall that a full
subcategory $\A$ of a category $\B$ is said to be {\it replete}
when every object of $\mathcal B$ isomorphic to an object  of
$\mathcal A$ is in  $\mathcal A$.)
\end{nada}

Now we can establish a basic fact of our theory which includes a
formal version of the Whitehead theorem in the homotopy theory of
topological spaces, and which is an easy corollary of Theorem
\ref{lemallave}. This theorem is no longer true with $\M[\mathcal
S^{-1}]$ in the place of  $\M[\mathcal S^{-1},\C]$ (see Remark
\ref{minimal no localiza}).

\begin{theorem}\label{W}  Let  $(\C,\mathcal S,\W) $ be a
category with strong and weak equivalences and  $\M$ be a full
subcategory of $\C_{cof}$. The functor $j$ induces a full and
faithful functor
$$\M[\mathcal{S}^{-1}, \C]\lra \C[\W^{-1}].$$ In particular this induced
functor reflects isomorphisms, that is to say, if
$w\in\C[\mathcal{S}^{-1}](M, N)$ is an isomorphism in
$\C[\W^{-1}]$, where $M$ and $N$ are in $\M$, then $w$ is an
isomorphism in $\loc{C}{S}$. \hfill$\Box$
\end{theorem}

\subsection{Cartan-Eilenberg categories}

For a category $\C$ with strong and weak equivalences the general
problem is to know if there are enough cofibrant objects. This
problem is equivalent to the orthogonal category problem for
$(\C[\mathcal S^{-1}],\delta(\W))$ (see \cite{Bo}(I.5.4)), which
has been studied by Casacuberta and Chorny  in the context of
homotopy theory (see \cite{CCh}).

\begin{definition}\label{defCE}
A category with strong and weak equivalences $(\C,\mathcal S,\W)$
is called a \emph{left Cartan-Eilenberg category} if each object
of $\C$ has a cofibrant
 left model (see Definitions \ref{defcofobj} and \ref{defmodel2}).

A category with weak equivalences $(\C,\W)$ is called a left
Cartan-Eilenberg category when the triple $(\C, \mathcal S,\W)$,
with $\mathcal S$ the class of isomorphisms of $\C$, is a left
Cartan-Eilenberg category.
\end{definition}

\begin{theorem}\label{CEimplies models}
 A  category with strong and weak
equivalences $(\C,\mathcal S,\W) $  is a left Cartan-Eilenberg
category if and only if
 $$j:\C_{cof}[\mathcal{S}^{-1},\C]\lra \loc{C}{W}$$
is an equivalence of categories.
\end{theorem}
\begin{proof}
By Theorem \ref{W}, $j$ is fully faithful. If $\C$ is a left
Cartan-Eilenberg, for each object $X$ there exists a cofibrant
left model  $\varepsilon: M\lra X$ of $X$, hence
$\gamma'(\varepsilon):M\lra X$ is an isomorphism in $\C[\W^{-1}]$,
so $j$ is essentially surjective.

Conversely, if $j$ is   an essentially surjective functor,  for
each object $X$, there exists a cofibrant object $M$ and an
isomorphism $\rho:M\lra X$ in $\C[\mathcal W^{-1}]$. By Theorem
\ref{lemallave}, there exists a morphism $\sigma:M\lra X$ in
$\C[\mathcal S^{-1}]$ such that $\gamma'(\sigma)=\rho$, therefore
$\sigma:M\lra X$ is a cofibrant left model of $X$, hence
$(\C,\mathcal S,\W) $ is a left Cartan-Eilenberg category.
\end{proof}

In a left Cartan-Eilenberg category the cofibrant left model is
functorial in the localised category $\C[\mathcal S^{-1}]$. More
precisely we have the following result.

\begin{corollary}\label{functor r} Let  $(\C,\mathcal S,\W) $ be a left Cartan-Eilenberg category.
There exists a functor
$$r:\loc{C}{S}\lra
\C_{cof}[\mathcal{S}^{-1},\C]$$ and a natural transformation
$$
\varepsilon':ir\Rightarrow 1
$$
such that:
\begin{enumerate}[(1)]
  \item For each object $X$, $\varepsilon'_X:ir(X)\lra X$ is a cofibrant left model of
  $X$.
  \item $r$ sends morphisms in $\delta(\W)$ into isomorphisms, and
  induces an equivalence of categories
 $$ \overline{r}: \loc{C}{W}\lra
\C_{cof}[\mathcal{S}^{-1},\C]$$quasi-inverse of $j$,  such that
$\overline{r}\gamma'=r$.

\item There exists a natural isomorphism
$\overline{\varepsilon}:j\overline{r}\Rightarrow 1_{\C[\W^{-1}]}$
such that $\gamma'\varepsilon'=\overline{\varepsilon}\gamma'$.

\item  The natural transformations
$$\gamma'\varepsilon':\gamma'ir\Rightarrow \gamma' ,\quad
\varepsilon'i:iri\Rightarrow i,\quad r\varepsilon':rir\Rightarrow
r$$ are isomorphisms.
\end{enumerate}

\end{corollary}
\begin{proof}
By the previous theorem, there exists a functor

$$
\overline{r}:
\loc{C}{W}\lra \C_{cof}[\mathcal{S}^{-1},\C]$$

that is the quasi-inverse of $j$, together with an isomorphism
$\overline{\varepsilon}:j\overline{r}\Rightarrow 1_{\C[\W^{-1}]}$.
Let

$$
r:=\overline{r}\gamma':\loc{C}{S}\lra
\C_{cof}[\mathcal{S}^{-1},\C].
$$

For each object $X$ in
$\loc{C}{S}$, $ir(X)$ is a cofibrant object, and
$\overline{\varepsilon}_{\gamma'X}:\gamma' ir X\lra \gamma' X$ is
an isomorphism in $\C[\W^{-1}]$, hence, by Theorem
\ref{lemallave}, there exists a unique morphism
$\varepsilon'_X:ir(X)\lra X$ in $\C[\mathcal S^{-1}]$ such that
$\gamma'(\varepsilon'_X)=\overline{\varepsilon}_{\gamma'X}$. If
$f:X\lra Y$ is a morphism in $\C[\mathcal S^{-1}]$, since
 $\overline{\varepsilon}$ is a natural
transformation, we have
$$\gamma'(f\circ \varepsilon'_X)
=\gamma'(f)\circ \overline{\varepsilon}_{\gamma' X}
=\overline{\varepsilon}_{\gamma' Y}\circ  \gamma'
i{r}(f)=\gamma'(\varepsilon'_Y\circ (ir)(f) ),$$ hence $f\circ
\varepsilon'_X= \varepsilon'_Y \circ (ir)(f)$, because $ir(X)$ is
cofibrant. As a consequence $\varepsilon':ir\Rightarrow 1$ is a
natural transformation. Therefore $\varepsilon'_X:ir(X)\lra X$ is
a functorial cofibrant left model of $X$.

 On the other
hand,$\gamma'\varepsilon'=\overline{\varepsilon}\gamma'$ and
$r\varepsilon'=\overline{r}\gamma'\varepsilon'=\overline{r}\overline{\varepsilon}\gamma'$
are  isomorphisms, since $\overline{\varepsilon}$ is an
isomorphism. By Theorem \ref{W},  $\varepsilon'i$ is also an
isomorphism.
\end{proof}

When proving that a category with strong and weak equivalences is
a Cartan-Eilenberg category, recognising cofibrant objects may
prove difficult, as the definition is given in terms of a lifting
property in $\C\mathcal [\mathcal S^{-1}]$.
 The sufficient
conditions we state in the next result are basic properties of the
category of bounded below chain complexes of modules over a
commutative ring in the Cartan-Eilenberg approach to homological
algebra (\cite{CE}).

These conditions are also the basic properties of the category of
 $\mathbf{k}$-cdg algebras in
Sullivan's theory of minimal models (see \cite{GM}). We followed
the same approach  to study the homotopy theory of modular operads
in  \cite{GNPR1}: see Theorem \ref{MOP} in this paper. We will also
apply it to  study the homotopy theory of filtered complexes
(see Theorem \ref{CEfiltrados}).

\begin{theorem}\label{F}
Let $(\C,\mathcal S,\W) $ be a category with strong and weak
equivalences  and $\M$ a full subcategory of $\C$. Suppose that

\begin{enumerate}[(i)]
\item for any $w:Y\longrightarrow X\in \W$ and any
$f\in\C(M, X)$, where $M\in Ob\,\M$, there exists a morphism
$g\in\loc{C}{S}(M, X)$ such that $w\circ g=f$ in $\loc{C}{S}$;

\item for any
$w:Y\longrightarrow X\in \W$ and any $M\in Ob\,\M$, the map
$$w_*:\loc{C}{S}(M,Y)\longrightarrow \loc{C}{S}(M,X)$$ is injective;
 and

\item for each object $X$ of $\C$ there exists a morphism
$\varepsilon:M\lra X$ in $\C$ such that $\varepsilon\in {\W}$
 and $M\in Ob\,\M$;

\end{enumerate}

Then,
\begin{enumerate}[(1)]
\item every object in $\M$ is cofibrant;

\item  $(\C,\mathcal S,\W)$ is a left Cartan-Eilenberg category;
and
\item the functor $\M[\mathcal{S}^{-1},
\C] \longrightarrow \C\left[\W^{-1}\right]$ is an equivalence of
categories.
\end{enumerate}

\end{theorem}
\begin{proof} Property  (2) follows  immediately from (1) and  (iii).
Property (3) follows from (iii), (1) and Theorem \ref{W}. So it is
enough to prove (1), that is: given $w:Y\lra X\in \W$,  $M$ in
$\M$ and $f\in \loc{C}{S}(M,X)$, there exists a unique $g\in
\loc{C}{S}(M,Y)$ such that $wg=f$ in $ \loc{C}{S}$. By (ii) it is
enough to prove the existence of $g$.

Suppose that $f\in \loc{C}{S}(M,X)$ can be represented as  an
alternating $\mathcal S$-zigzag of $\C$ of  length $m$, from $M$
to $X$. We proceed by induction on $m$. The case $m=1$  follows
from hypothesis (i).

Let $m>1$.  Then  $f=f_2s^{-1}f_1$, where $f_1\in \C(M,X_1)$,
$s:X_2\longrightarrow X_1\in \mathcal S$ and
$f_2:X_2\longrightarrow X$ is an alternating $\mathcal S$-zigzag
of $\C$ of length $m-2$. By (iii), there exists a morphism
$\varepsilon:M_2\longrightarrow X_2$ in ${\W}$ such that $M_2\in
Ob\,\M$, hence, by
 (i),  there
exists $g_1\in \loc{C}{S}(M,M_2)$ such that $f_1=s\varepsilon
g_1$. In addition, by the induction hypothesis, since
$f_2\varepsilon$ can be represented as an alternating
 $\mathcal S$-zigzag of $\C$ of length $m-2$,  there exists $g_2\in
\loc{C}{S}(M_2,Y)$ such that $f_2\varepsilon=wg_2$. Then
$g:=g_2g_1\in \C[\mathcal S^{-1}](M,Y)$ satisfies  $wg=f$.
$$
\xymatrix{
& &M_2 \ar[d]^\varepsilon\ar@{.>}[r]^{g_2}&Y\ar[d]^w \\
  M \ar@{.>}[urr]^{g_1} \ar[r]^{f_1} &
  X_1&X_2\ar[l]_{s}\ar[r]^{f_2}&X}
$$
\end{proof}

\begin{example}\label{EjemploCE}
Let $\A$ be an abelian category with enough projective objects and
let $\mathbf C_+(\A)$ be the category of bounded below chain
complexes of $\A$. Let $\mathcal S$ be the class of homotopy
equivalences, and $\W$ the class of quasi-isomorphisms. Let $\M$
be the full subcategory of projective degree-wise complexes.
Because the localisation $\mathbf C_+(\A)[\mathcal S^{-1}]$ is the
homotopy category $\mathbf K_+(\A)$, by Proposition
\ref{congruencia} and Example \ref{cilindro}, the hypothesis of
the previous theorem are well known facts (see \cite{CE} and
\cite{GMa}), hence $(\mathbf C_+(\A),\mathcal S,\W)$ is a left
Cartan-Eilenberg category and $\M$ is a subcategory of cofibrant
left models of $\mathbf C_+(\A)$.

\end{example}

\subsection{Idempotent functors and
reflective subcategories} In some cases, localisation of
categories may be realised through reflective subcategories or,
equivalently, by Adams idempotent functors (see \cite{Bo}(3.5.2)
and  \cite{A}, section 2). These notions are also related with the
Bousfield localisation (see \cite{N} for this notion in the
context of triangulated categories).  The following Theorem
\ref{sufcof} relates left Cartan-Eilenberg categories with the
dual notions of coreflective subcategories and coidempotent
functors. Some of the parts of the theorem are a reinterpretation
of well known results when $\mathcal S$ is the trivial class of
the isomorphisms, which is in fact the key of the problem. For
triangulated categories, the fourth condition in Theorem
\ref{sufcof} corresponds to the notion of Bousfield colocalization
(see \cite{N}).

>From now on, we will use also the notation $*$ for the Godement
product between  natural transformations and functors (see
\cite{G}, Appendice), and apply its properties freely.

We recall that a replete subcategory (see \ref{+notacions})
 $\A$ of a category $\B$ is
called \emph{coreflective} if the inclusion functor $i:\A\lra \B$
admits a right adjoint $r:\B\lra\A$,  called a \emph{
coreflector}. We recall also that  a \emph{ coidempotent} functor
on a category $\B$ is a pair $(R,\varepsilon)$, where $R:\B\lra
\B$ is an endofunctor of  $\B$ and $\varepsilon$ is a morphism
$\varepsilon: R\Rightarrow 1_\B$, called \emph{counit}, such that
$$R\varepsilon, \varepsilon R:R^2\Rightarrow R$$ are  isomorphisms,
and $R\varepsilon=\varepsilon R$ (see \cite{A}). In fact, the
equality $R\varepsilon=\varepsilon R$ is a consequence of the
first condition, as proved in the following lemma.

\begin{lemma}\label{hrf}
Let $\B$ be a category together with an  endofunctor $R:\B\lra
\B$ and a  morphism $\varepsilon: R\Rightarrow 1_\B$ such that
the morphisms
$$\varepsilon R,\, R \varepsilon:R^2\Rightarrow R$$ are isomorphisms.
Then
  $(R,\varepsilon)$ is a coidempotent
  functor on $\B$.
\end{lemma}
\begin{proof}
In the  strict simplicial object associated to $(R,\varepsilon)$
(see \cite{G}, App.),
$$\xymatrix{\cdots
R^3\ar@<0.5ex>[r]\ar[r]\ar@<-0.5ex>[r]&R^2\ar@<0.5ex>[r]\ar@<-0.5ex>[r]&R,
}$$ with face morphisms
$$
\delta^n_i=R^i*\varepsilon*R^{n+1-i}:R^{n+1}\lra R^{n}, \quad 0\le
i\le n, \; 1\le n,
$$
the arrows $\delta_0^1=\varepsilon R,\delta_1^1=R \varepsilon $
are isomorphisms. From the simplicial relations
$$
\delta_0^1\delta_0^2=\delta_0^1\delta_1^2,\;
\delta_1^1\delta_1^2=\delta_1^1\delta_2^2
$$
we deduce $\delta_0^2=\delta_1^2=\delta_2^2$. Since
$\delta_0^1\delta_2^2=\delta_1^1\delta_0^2,\;$ and
$\delta_2^2=\delta^2_0=\varepsilon R^2$ is also an isomorphism, we
conclude that $\delta_0^1=\delta_1^1$.
\end{proof}

\begin{theorem}\label{sufcof}
Let $(\C,\mathcal S,\W)$ be a category with strong and weak
equivalences. Then the following conditions are equivalent.

\begin{itemize}
\item[(i)]  $(\C, \mathcal S,\W)$ is a left Cartan-Eilenberg
category.

\item[(ii)]
There exists a coidempotent functor $(R',\varepsilon')$ on
$\loc{C}{S}$ such that $\overline{\W}$ is the pre-image by $R'
\delta$ of the class of isomorphisms in $\C[\mathcal S^{-1}]$, and
$\gamma'\varepsilon'$ is an isomorphism.

\item[(iii)]  The inclusion functor $i:\C_{cof}[\mathcal{S}^{-1},\C]
\lra\loc{C}{S}$ admits a right adjoint
$$r:\loc{C}{S}\lra\C_{cof}[\mathcal{S}^{-1},\C],$$ with a counit
$\varepsilon':i r\Rightarrow 1$,  such that
$\overline{\delta(\W)}$ is the pre-image by $r$ of the class of
isomorphisms in $\C_{cof}[\mathcal{S}^{-1},\C]$, and
$r\varepsilon'$ is an isomorphism. In particular
$\C_{cof}[\mathcal{S}^{-1},\C]$ is a coreflective subcategory of
$\loc{C}{S}$.

\item[(iv)]
The localisation functor $\gamma':\C[\mathcal S^{-1}]\lra
\C[\W^{-1}]$ admits a left adjoint $$ \lambda:\C[\W^{-1}]\lra
\loc{C}{S}.$$
\end{itemize}
Assuming  that these conditions are satisfied,
$\C_{cof}[\mathcal{S}^{-1},\C]$ is the essential image of $R'$
(and $\lambda$).

\end{theorem}
\begin{proof}
We prove the theorem in several steps. Firstly we recall, from
Corollary \ref{functor r}, that if $(\C, \mathcal S,\W)$ is a left
Cartan-Eilenberg category there exists a  functor
$$r:\loc{C}{S}\lra\C_{cof}[\mathcal{S}^{-1},\C],$$ together with a
morphism $\varepsilon':i r \Rightarrow 1$ such that
$\varepsilon'*i$, $r*\varepsilon'$ and $\gamma'*\varepsilon'$ are
isomorphisms.

Step $1$: (i) implies (ii).
 Let $R':\C[\mathcal S^{-1}]\lra \C[\mathcal S^{-1}]$
be the functor $R'=ir$. Then $\varepsilon':R'\Rightarrow 1$ is a
natural transformation, and $\varepsilon'*R'=\varepsilon'*(i
r)=(\varepsilon'*i)*r$ and $R'*\varepsilon'=(i
r)*\varepsilon'=i*(r*\varepsilon')$ are
 isomorphisms, because so are $\varepsilon'*i$ and
$r*\varepsilon'$. Therefore, by Lemma \ref{hrf},
$(R',\varepsilon')$ is a coidempotent functor.

Let us see that  $\overline{\W}$ is the pre-image by $R'\delta$ of
the class of isomorphisms in $\C[\mathcal S^{-1}]$. It is enough
to see that, given a morphism $f:X\lra Y$ in $\C[\mathcal
S^{-1}]$, $R'(f)$ is an isomorphism if and only if $\gamma'(f)$ is
an isomorphism. From the naturality of $\varepsilon'$ we have
$$
\varepsilon'_Y\circ R'(f)=f\circ\varepsilon'_X \ ,
$$
therefore,  by  Theorem $\ref{W}$, $\gamma'(f)$ is an isomorphism
if and only if $R'(f)$  is an isomorphism.

 Step $2$: (i)  implies (iii).   For each category $\mathcal X$, the functor
 $$i_*: \Cat(\mathcal X,\C_{cof}[\mathcal S^{-1},\C])
 \lra \Cat(\mathcal X,\C[\mathcal S^{-1}]) $$ is fully faithful;
hence, to define a natural transformation \; $\eta: 1\Rightarrow r
i$, it is enough to define a natural transformation\;
$i*\eta:i\Rightarrow i ri$. Since $\varepsilon'*i:iri\Rightarrow
i$ is an isomorphism, we define $\eta$ to be  such that
$i*\eta=\left(\varepsilon'*i\right)^{-1}$. Let us check that
$\eta$ and $ \varepsilon'$ are the unit and  the counit,
respectively, of an adjunction $i\vdash r$, that is to say (see
for example \cite{ML}),
$$
(r*\varepsilon')\circ (\eta*r)=1_r,\quad (\varepsilon'*i)\circ
(i*\eta)=1_i.
$$
By step 1, $(ir)*\varepsilon'=\varepsilon'*(ir)$, and by the
definition of $\eta$ we obtain
$$\begin{aligned}
  &i*((r*\varepsilon')\circ (\eta*r))=((ir)*\varepsilon')\circ
(i*\eta*r) =(\varepsilon'*(ir))\circ(
\left(\varepsilon'*i\right)^{-1}*r) \\
=& ((\varepsilon'*i)*r)\circ \left((\varepsilon'*i\right)^{-1}*r
)= \left((\varepsilon'*i)\circ (\varepsilon'*i)^{-1}\right)*r=
1_{i}*r=i*1_r.
\end{aligned}
$$
Since $i_*$ is fully faithful, we obtain $(r*\varepsilon')\circ
(\eta*r)=1_r$. The other identity being trivial, we conclude that
$r$ is a right adjoint for $i$.

The other assertions are  consequence of step $1$.

 Step 3: (i) implies (iv). By Corollary \ref{functor r} there is a
functor $\overline{r}:\C[\W^{-1}]\lra \C_{cof}[\mathcal S^{-1},\C]
$ such that $\overline{r} \gamma'= r$. Let $\lambda=i
\overline{r}$. Since
$$
\gamma'^*:\Cat(\C[\W^{-1}],\C[\W^{-1}])\lra \Cat(\C[\mathcal
S^{-1}],\C[\mathcal W^{-1}])
$$
is fully faithful, and $\gamma'*\varepsilon':
\gamma'\lambda\gamma'\Rightarrow \gamma'$ is an isomorphism, there
exists a unique morphism $\eta:1\Rightarrow \gamma' \lambda$ such
that
$$
\eta*\gamma'=(\gamma'*\varepsilon')^{-1}.
$$
Then, $(\eta, \varepsilon')$ are the unit and the counit of an
adjunction $\lambda\dashv \gamma'$, that is to say,
$$
(\gamma'*\varepsilon')\circ (\eta*\gamma')=1_{\gamma'},\quad
(\varepsilon'*\lambda)\circ (\lambda*\eta)=1_{\lambda}.
$$
Indeed, the first identity follows trivially from the definition
of $\eta$. For the second one, we have $\lambda\gamma'=i
\overline{r} \gamma'=i r$ by the definitions, and
$(ir)*\varepsilon'=\varepsilon'*(ir)$ by step $1$, so we have
$$\begin{aligned}
&\left((\varepsilon'*\lambda)\circ
(\lambda*\eta)\right)*\gamma'=(\varepsilon'*(\lambda
\gamma'))\circ (\lambda*\eta*\gamma')
  =(\varepsilon'*(ir))\circ ((i
\overline{r})*(\gamma'*\varepsilon')^{-1})\\=&
((ir)*\varepsilon')\circ ((i
\overline{r})*(\gamma'*\varepsilon')^{-1})
=((i\overline{r})*(\gamma'*\varepsilon'))\circ
((i\overline{r})*(\gamma'*\varepsilon')^{-1})
\\=&(i\overline{r})*((\gamma'*\varepsilon')\circ
(\gamma'*\varepsilon')^{-1})=
(i\overline{r})*1_{\gamma'}=\lambda*1_{\gamma'}=1_{\lambda}*\gamma',
\end{aligned}
$$
therefore, since $\gamma'^*$ is fully faithful,  the second
identity of the adjunction is also satisfied.

Step 4: (ii) implies (i). Firstly, for each object $X$, let us
check that $R'X$ is cofibrant.  Let
  $w:A\lra B$ be a morphism in  $\delta(\W)$. By hypothesis
$R'( w)$ is an isomorphism, therefore we have a commutative
diagram
$$
\xymatrix{ \C[\mathcal
S^{-1}](R'X,R'A)\ar[dd]^{R'w_*}\ar[rr]^{\varepsilon'_{A*}}&&
\C[\mathcal S^{-1}](R'X,A)\ar[dd]^{w_*}\\
\\\C[\mathcal S^{-1}](R'X,R'B)\ar[rr]^{\varepsilon'_{B*}}
&&\C[\mathcal S^{-1}](R'X,B)
 }
$$
where $\varepsilon'_{A*}$, $\varepsilon '_{B*}$ and $R'w_*$ are
bijective. Therefore  $w_*$ is bijective, thus $R'X$ is cofibrant.
Since $\varepsilon'_X:R'(X)\lra X \in \overline{\delta(\W)}$, each
object has a cofibrant left model, hence $(\C,\mathcal S,\W)$ is a
left Cartan-Eilenberg category.

Step 5: (iii) implies (i). For each object $X$,
$\varepsilon'_X:ir(X)\lra X$ is a cofibrant left model of $X$,
therefore $(\C,\mathcal S,\W)$ is a left Cartan-Eilenberg
category.

 Step 6: (iv) implies (i).
This is an easy consequence of the dual of Proposition I.1.3 of
\cite{GZ}. In fact, let  $\eta:1\Rightarrow \gamma'\lambda$ and
$\varepsilon':\lambda\gamma'\Rightarrow 1$ be the unit and the
counit of the adjunction, respectively. The  functor $\C[\mathcal
S^{-1}][\delta(\W)^{-1}]\lra \C[\W^{-1}]$ induced by $\gamma'$ is
an isomorphism, thus, by \emph{loc. cit.}, $\eta$ is an
isomorphism. Therefore the identity of the adjunction
$$
(\eta*\gamma')\circ (\gamma'*\varepsilon')=1_{\gamma'}
$$
proves that $\gamma'*\varepsilon'$ is an isomorphism. So, for each
object $X$, $\varepsilon'_X:\lambda\gamma'(X)\lra X$ is a left
model.  On the other hand, for each pair of objets $X$ and $Y$,
the composition
$$
\xymatrix{\C[\mathcal
S^{-1}](\lambda\gamma'(X),Y)\ar[r]^{\gamma'_Y}&
\C[\W^{-1}](\gamma'\lambda\gamma'(X),\gamma'(Y))
\ar[r]^{\eta_{\gamma'(X)}^*}&\C[\W^{-1}](\gamma'(X),\gamma'(Y))}
$$
is the adjunction map, and as $\eta_{\gamma'(X)}^*$ is bijective,
so is  $\gamma'_Y$. Therefore, by Proposition \ref{lemallave},
$\lambda\gamma'(X)$ is cofibrant. Hence,
$\varepsilon'_X:\lambda\gamma'(X)\lra X$ is a cofibrant left model
of $X$, which proves (i).

Finally, in step $5$ (resp. step $6$) we have just proved that
$R'X$ (resp. $\lambda\gamma'(X)$) is cofibrant, for each object
$X$. Conversely, if $M$ is cofibrant, $\varepsilon'_M:R'M\lra M$
(resp. $\varepsilon'_M:\lambda \gamma'M\lra M$) is a morphism in
$\overline{\delta(\W)}$ between cofibrant objects, therefore, by
Theorem \ref{W}, it is an isomorphism in $\C[\mathcal S^{-1}]$. So
$\C_{cof}[\mathcal S^{-1},\C]$ is the essential image of $R'$
(resp. $\lambda$).
\end{proof}

\begin{nada}\label{notacions}
Let $(\C,\mathcal S,\W)$ be a left Cartan-Eilenberg category. We
 summarise the different functors  we have
encountered  between the categories associated to $(\C,\mathcal
S,\W)$ in the following diagram

$$
\xymatrix{\C\ar[dd]^\gamma\ar[rr]^\delta&&\C[\mathcal
S^{-1}]\ar@<0.5ex>[ddll]^{\gamma'}\ar@<0.5ex>[dd]^r\\
\\
\C[\W^{-1}]\ar@<0.5ex>[uurr]^\lambda\ar@<-0.5ex>[rr]
_{\overline{r}}&&\C_{cof}[\mathcal
S^{-1},\C]\ar@<-0.5ex>[ll]_j\ar@<0.5ex>[uu]^i,}
$$
where: \begin{itemize} \item[(a)] The functors $\gamma$, $\delta$
and $\gamma'$ are the localisation functors (see
\ref{notations2}).

\item[(b)] The functor $i$ is the inclusion functor (see
\ref{+notacions}) and $r$ is the functorial cofibrant left model
(see Theorem \ref{functor r}).
\item[(c)] The functor $r$ is the right adjoint
of $i$ (see Theorem \ref{sufcof} (iii)).

\item[(d)]
The functor $\overline r$ is the unique functor such that
$r={\overline r}\gamma '$.

\item[(e)]
The functor $j$ is defined by    $j:=\gamma ' i$ (see
 \ref{+notacions}).

\item[(f)]
The functors $j$ and $\overline{r}$ are quasi-inverse equivalences
(see Corollary \ref{functor r}).

\item[(g)]
The functor $\lambda$ is  defined by $\lambda:=i \overline r$. It
is left adjoint for $\gamma'$ (see Theorem \ref{sufcof}, (iv)).
\end{itemize}
\end{nada}

\begin{remark}\label{horesolvent=resolvent}
If $\mathcal S$ is just the class of isomorphisms, then $\C_{cof}$
is the class of objects which are left orthogonal (see
\cite{Bo}(5.4))  to $\W$, therefore $(\C,\W)$ is a left
Cartan-Eilenberg category if and only if $\C_{cof}$ is a
coreflective subcategory of $\C$.

\end{remark}

\subsection{Resolvent functors}

Sometimes the coidempotent functor $R':\C[\mathcal S^{-1}]\lra
\C[\mathcal S^{-1}]$  in Theorem \ref{sufcof} comes from an
endofunctor of $\C$ itself. We formalise this situation in the
following definition.

\begin{definition}
Let  $(\C,\mathcal S,\W)$ be a category with strong and weak
equivalences. A \emph{left resolvent functor} on $\C$ is a pair
$(R,\varepsilon)$ where
\begin{enumerate}[(i)]
\item   $R:\C\longrightarrow \C$ is a functor such that $R(X)$ is
a cofibrant object, for each $X\in Ob\; \C$; and

\item  $\varepsilon:R\Rightarrow \id_\C$ is morphism such that
$\varepsilon_X:R(X)\longrightarrow X$ is in $\overline{\W}$, for
each $X\in Ob\, \C$. \end{enumerate}\end{definition}

A left resolvent functor is also called a functorial cofibrant
replacement.

\begin{lemma}\label{LFR} Let  $(\C,\mathcal S,\W)$ be a category with
strong and weak equivalences, and let $(R,\varepsilon)$ be a left
resolvent functor on $\C$. Then,

\begin{enumerate}

\item we have $\overline{\W}=R^{-1}(\overline {\mathcal S})$, in
particular $R(\mathcal S)\subset \overline{\mathcal S} $;

\item  we have
$R(\varepsilon_X),\varepsilon_{R(X)}\in \overline{\mathcal S}$,
for each  $X\in \text {Ob\;} \C$; and

\item $(R,\varepsilon)$ induces a coidempotent functor
$(R',\varepsilon')$ on $\C[\mathcal S^{-1}]$.

\end{enumerate}
\end{lemma}
\begin{proof}
Since $R^{-1}( \overline{\mathcal S})$ is a saturated class of
morphisms, in order to prove that $\overline{\W}\subset
R^{-1}(\overline{\mathcal S}) $ it is enough to check that
$\W\subset R^{-1}( \overline{\mathcal S})$. In fact, if $w:X\lra
Y$ is a morphism in $\W$, we have a commutative diagram
$$\xymatrix{R(X)\ar[rr]^{R(w)}\ar[dd]_{\varepsilon_X}&&
R(Y)\ar[dd]^{\varepsilon_Y}\\ \\X\ar[rr]^w&&\, Y\, , }
$$
where $w$, $\varepsilon_X$ and $\varepsilon_Y$ are morphisms in
$\overline{\W}$, hence $R(w)$ is also in $\overline{\W}$, since
$\overline{\W}$ has the \emph{2 out of 3} property. By Theorem
\ref{W}, $R(w)$ is in $\overline{\mathcal S}$, therefore
${\W}\subset R^{-1}( \overline{\mathcal S}) $. Conversely, if
$w\in R^{-1}(\overline{\mathcal {S}}) $, then $R(w)\in
\overline{\mathcal S}$, and, from the previous diagram, we obtain
$w\in \overline{\W}$.

>From the hypothesis and part (1) we obtain $R\varepsilon_X\in
\overline{\mathcal S}$. Next, from $\varepsilon_{RX}\in
\overline{\W}$ and Theorem \ref{W}, we obtain $\varepsilon_{RX}\in
\overline{\mathcal S}$. Finally (3) follows  from (2) and Lemma
\ref{hrf}.
\end{proof}

 A category with a left
resolvent functor  is a particular left Cartan-Eilenberg category
where both localisations $\C_{cof}[\mathcal S^{-1},\C]$ and
$\C_{cof}[\mathcal S^{-1}]$ agree.

\begin{proposition}\label{PreMA}
 Let  $(\C,\mathcal S,\W)$ be a category with strong and weak
equivalences, and let $(R,\varepsilon)$ be a left resolvent
functor on $\C$. Then,
\begin{enumerate}
\item  $(\C,\mathcal S,\W)$
is a left Cartan-Eilenberg category;

\item the canonical functor $\alpha:\C_{cof}[\mathcal S^{-1}]
\lra \C[\W^{-1}]$ is an
equivalence of categories; and
\item an object $X$ of $\C$ is cofibrant if and only if
$\varepsilon_X:RX\lra X$ is   an isomorphism in $\C[\mathcal
S^{-1}]$.

\end{enumerate}
\end{proposition}

\begin{proof} First
of all, for each object $X$ of  $\C$, we have
$\varepsilon_X:RX\lra X\in\overline{\W}$, where $RX$ is cofibrant.
In particular, $\varepsilon_X:RX\lra X$ is a cofibrant left model
of $X$, therefore $\C$ is a left Cartan-Eilenberg category, which
proves (1).

Next, let us see (2). Since $R(X)$ is cofibrant and $R(\W)\subset
\overline{\mathcal S}$, by Lemma \ref{LFR}, the functor $R$
induces a functor $$\beta:\C[\W^{-1}]\lra \C_{cof}[\mathcal
S^{-1}]$$ such that $\delta R=\beta\gamma$. Let us see that
$\beta$ is a quasi-inverse of $\alpha$. Indeed, for each object
$X$ of $\C$, the counit $\varepsilon_X:R(X)\lra X$ induces a
morphism in $\C[\W^{-1}]$
$$\gamma(\varepsilon_X):\alpha\beta(\gamma(X))=\gamma(R(X))\lra
\gamma(X)$$ which is an isomorphism. On the other hand, for each
cofibrant object $M$, the morphism
$$\delta(\varepsilon_M):\beta\alpha(\delta(M))=\delta(R(M))\lra
\delta(M)$$ satisfies
$\alpha\delta(\varepsilon_M)=\gamma(\varepsilon_M)$, so it is an
isomorphism. Therefore, by    Theorem \ref{W}, $\varepsilon_M\in
\overline{\mathcal S}$. So $\delta(\varepsilon_M)$ is an
isomorphism, which proves (2).

Finally, since $R$ is a left resolvent functor, $R(X)$ is a
cofibrant object for each object $X$, hence, if $\varepsilon_X$ is
an isomorphism in $\C[\mathcal S^{-1}]$, $X$ is also cofibrant.
Conversely, if $X$ is cofibrant, then $\varepsilon_X:RX\lra X$ is
a morphism in $\overline{\W}$
 between cofibrant objects, hence, by
Theorem \ref{W}, it is  an isomorphism in $\C[\mathcal S^{-1}]$.
\end{proof}

 The following result gives a useful criterion in
order to obtain left resolvent functors, as we will see in section
$6$.

\begin{theorem}\label{MA}
Let $\C$ be a category, $\mathcal S$ a class of morphisms in $\C$,
$R:\C\longrightarrow \C$ a functor and $\varepsilon:R\Rightarrow
\id$ a morphism such that
$$
R(\mathcal S)\subset \mathcal S,\quad R(\varepsilon_X)\in \mathcal
S, \quad \varepsilon_{R(X)}\in \mathcal S,
$$
for each $X\in \text {Ob\;} \C$. If we take
$\W=R^{-1}(\overline{\mathcal S})$, then $\mathcal S\subset \W$
and $(R,\varepsilon)$ is a left resolvent functor for
$(\C,\mathcal S,\W)$,  which is therefore a left Cartan-Eilenberg
category satisfying conditions (1), (2) and (3) of Proposition
\ref{PreMA}.

\end{theorem}

\begin{proof} The pair  $(R,\varepsilon)$ induces a coidempotent functor
$(R',\varepsilon')$ on $\loc{C}{S}$ which satisfies the hypothesis
(ii) of Theorem \ref{sufcof},
 therefore $\varepsilon_X:R(X)\lra X$ provides a
cofibrant left model of $X$, for each $X$. Hence $(R,\varepsilon)$
is a left resolvent functor for $(\C,\mathcal S,\W)$.
\end{proof}

\begin{example}\label{ComplejosSonCE}
Let  $\mathbf C_+(A)$ be the category of bounded below chain
complexes of $A$-modules, where $A$ is a commutative ring, and let
$\mathcal S$ be the class of homotopy equivalences. Let $R$ be the
endofunctor on $\mathbf C_+(A)$ defined by the free functorial
resolution induced by the functor on the category of $A$-modules,
$X\mapsto  A^{(X)}$, where $A^{(X)}$ denotes the free $A$-module
with base $X$, and $\varepsilon:R\Rightarrow \id$ is the
augmentation morphism. Since the objects of $\mathbf C_+(A)$ are
bounded below chain complexes, a quasi-isomorphism between two such
complexes which are free component-wise is a homotopical
equivalence. Hence the hypothesis of the previous theorem are
verified and, therefore, $(R,\varepsilon)$ is a left resolvent
functor on $\mathbf C_+(A)$. Moreover, the class $\W$ is the class
of quasi-isomorphisms (as in Example \ref{EjemploCE}), and the
cofibrant objects are the complexes which are homotopically
equivalent to  a free component-wise complex.

In the next sections 3 and 6 we will see other examples of
resolvent functors.
\end{example}

\begin{remark}
The dual  notions  of cofibrant object and left Cartan-Eilenberg
category, are the notions of \emph{fibrant object} and \emph{right
Cartan-Eilenberg category}. All the preceding results have their
corresponding dual. For example, dual of Theorem \ref{CEimplies
models} says that a category with strong and weak equivalences
$(\C,\mathcal S,\W)$ is a right Cartan-Eilenberg category if and
only if  the  functor $\C_{fib}[\mathcal S^{-1},\C] \lra
\loc{C}{W}$ is an equivalence of categories.
\end{remark}

\section{Models of functors and derived  functors}

In this section we study functors defined on a Cartan-Eilenberg
category $\C$ and taking values in a category $\D$ with a class of
weak equivalences. We prove that, subject to some hypotheses,
certain categories of functors are also Cartan-Eilenberg
categories.  In this context we can realise derived functors, when
they exist, as cofibrant models in the functor category. The
classic example is the category of additive functors defined on
a category of complexes of an abelian category with enough
projective objects.

\subsection{Derived functors}
To begin with, we recall the definition of derived functor as set
up by Quillen (\cite{Q}).

 Let  $(\C,\W)$ be a category with weak
equivalences, and $\D$ an arbitrary  category. Recall that the
category $\Cat(\C[\W^{-1}],\D)$ is identified, by means of the
functor
$$
\gamma^*:\Cat(\C[\W^{-1}],\D)\lra \Cat(\C,\D),
$$
with the full subcategory $\Cat_\W(\C,\D)$ of $\Cat(\C,\D)$ whose
objects are the functors which send morphisms in $\W$ to
isomorphisms in $\D$.

If $F:\C\lra \D$ is a functor, a right Kan extension (see
\cite{ML}, Chap. X) of $F$ along $\gamma:\C\lra \C[\W^{-1}]$ is a
functor
$$\Ran_\gamma F:\C[\W^{-1}]\lra \D,$$ together with a natural
transformation $\theta_F=\theta_{\gamma,F}:(\Ran_\gamma
F)\gamma\Rightarrow F$, satisfying the usual universal property.

\begin{definition}\label{DefL}
Let  $(\C,\W)$ be a category with weak equivalences, and $\D$ an
arbitrary  category. A functor $F:\C\lra \D$ is called \emph{left
derivable} if it exists the right Kan extension of $F$ along
$\gamma$. The functor
$$\mathbb L_\W F:=(\Ran_\gamma F) \gamma$$  is called a \emph{left
derived functor of $F$} with respect to $\W$.

We will denote by  $\Cat'((\C,\W),\D)$  the full subcategory of
$\Cat(\C,\D)$ of  left derivable functors with respect to $\W$.
\end{definition}
\begin{nada}

The left derived functor $\mathbb L_\W F$  is endowed with a
natural transformation $\theta_F:\mathbb L_\W F\Rightarrow F$ such
that, for each functor $G\in Ob\Cat_\W(\C,\D)$ the map
$$
\mathbf {Nat}(G,\mathbb L_\W F)\lra \mathbf{Nat}(G,F),\quad
\phi\mapsto \theta_F\circ \phi
$$
is bijective.
\end{nada}

If $\W$ has a right calculus of fractions, the definition of left
derived functor agrees with the definition given by Deligne in
\cite{D2}.

Functors in $\Cat_\W(\C,\D)$ are tautologically derivable functors
as ensues from the following easy lemma.

\begin{lemma}\label{lema1} Let  $(\C,\W)$ be a category
with weak equivalences, and $\D$ an arbitrary category. Then,
\begin{enumerate}

\item  any functor $F : \C \longrightarrow \D$ which takes
$\mathcal{W}$ into isomorphisms  induces a unique functor
$$F':\C[\W^{-1}]\lra\D$$ such that $F'\gamma=F$. This functor
$F'$ satisfies  $F'=Ran_\gamma F$, with  $\theta_F=Id$. In
particular, $F$   is left derivable and $\mathbb L_\W F=F$; and

\item  $\Cat_\W(\C,\D)$ is a full subcategory of
$\Cat'((\C,\W),\D)$. \hfill$\Box$

\end{enumerate}

\end{lemma}
\begin{nada}
For each $F\in  Ob\,\Cat'((\C, \W),\D)$, we have $\mathbb L_\W
F\in Ob\,\Cat_\W(\C,\D)$, so, by the previous lemma, part (1), it
results that $\mathbb L_\W F\in Ob\, \Cat'((\C,\W),\D)$.
Therefore, taking the left derived functor $\mathbb L_{\W}$
defines a functor $$\mathbb L_{\W} :\Cat'((\C,\W),\D)\lra
\Cat'((\C,\W),\D),$$ and the canonical morphism $\theta_F:\mathbb
L_{\W} F\lra F$ gives a natural transformation $\theta:\mathbb
L_{\W}\Rightarrow \id$.
\end{nada}

\begin{theorem}\label{derivargenerico} With the  notation above we have

(1)  the pair $(\mathbb L_\W,\theta)$ is a coidempotent functor on
$\Cat'((\C,\W),\D)$;

(2) the category with weak equivalences
$(\Cat'((\C,\W),\D),\widetilde{\W})$, where $\widetilde{\W}$ is
the class of morphisms whose image by $\mathbb L_\W$ is an
isomorphism, is a left Cartan-Eilenberg category; and

(3) the category $\Cat_\W(\C,\D)$ is the subcategory of its
cofibrant objects.

In particular, if $F:\C\lra \D$ is a left derivable functor, a
left derived functor of $F$ is the same as a cofibrant left model
of $F$.
\end{theorem}

\begin{proof} In the sequel we shorten $\mathbb L_\W$ as $\mathbb L$.
First of all, by Lemma \ref{lema1}, for each left derivable
functor $F:\C\lra \D$,  $\mathbb L \mathbb L F=\mathbb L F$ and
$\theta_{\mathbb L F}$ is the identity, hence $\theta_{\mathbb L
F}$ is an isomorphism. On the other hand, the naturality of
$\theta$ implies that the following diagram is commutative
$$
\xymatrix{ \mathbb L^2F\ar[rr]^{\theta_{\mathbb
LF}}\ar[dd]_{\mathbb L\theta_{F}}
&&\mathbb LF\ar[dd]^{\theta_F}\\
\\
\mathbb LF\ar[rr]^{\theta _F}\ar[rr]&&F\,,}
$$
hence, by the universal property of Definition \ref{DefL}, we
obtain $\mathbb L (\theta_F)=\theta_{\mathbb L F}$, so $\mathbb L
(\theta_F)$ is also an isomorphism. Therefore  $(\mathbb
L,\theta)$ is a coidempotent functor on $\Cat'((\C,\W),\D)$. So,
by Theorem \ref{sufcof}, $\Cat'((\C,\W),\D)$, is a left
Cartan-Eilenberg category, taking the isomorphisms as strong
equivalences, and  the class of morphisms of $\Cat'((\C,\W),\D)$
whose image by $\mathbb L$ is an isomorphism as weak equivalences.
Finally, the cofibrant objects are the functors isomorphic to
functors $\mathbb L F$, that is to say, the functors in
$\Cat_\W(\C,\D)$.
\end{proof}

\subsection{A derivability criterion for functors}
In this section we give  a derivability criterion for functors
defined on a left Cartan-Eilenberg category, which is a non
additive extension of the standard derivability criterion for
additive functors, and  we obtain a Cartan-Eilenberg category
structure for functors satisfying such derivability criterion.

In the following results we use the notation settled in
\ref{notacions}.

\begin{theorem}\label{derivarS}
Let $(\C,\mathcal S,\W)$ be a left Cartan-Eilenberg category. For
any category $\D$,
\begin{enumerate}
\item
$\Cat_{\mathcal S}(\C,\D)$ is a full subcategory of
$\Cat'((\C,\W),\D)$;

\item if $F\in Ob\ \Cat_{\mathcal S}(\C,\D)$, then
$$
\mathbb L_\W F =F' \lambda \gamma,
$$
where $F':\loc{C}{S}\lra \D$ denotes the functor induced by $F$;
and the canonical morphism $\theta_F:\mathbb L_\W F\lra F$ is
defined by $\theta_F=F'*\varepsilon*\delta$, that is to say,
$$
(\theta_F)_X=F'(\varepsilon'_{\delta X}),
$$
for each object $X$ of $\C$.

\end{enumerate}

\end{theorem}
\begin{proof}
By (iv) of  Theorem \ref{sufcof}, $\lambda$ is left adjoint to
$\gamma '$, and  $\varepsilon':\lambda \gamma'\Rightarrow 1$ is
the counit of the adjunction, therefore these functors induce a
pair of functors
$$
\xymatrix{ {\Cat(\C[\mathcal
S^{-1}],\D)}\ar@<-0.5ex>[rr]_{\lambda^\ast}
    & & \Cat(\C[\W^{-1}],\D), \ar@<-0.5ex>[ll]_{\gamma
'^\ast} }
$$
which are also adjoint, where $\lambda^*$ is right adjoint to
$\gamma^{'*}$, and $\varepsilon'^{*}:
\gamma'^*\lambda^*\Rightarrow 1$ is the counit of the adjunction,
as is easily seen. Hence, for each functor $G\in \Cat(\C[\mathcal
S^{-1}],\D)$, $\lambda^*(G)=G \lambda$ is a right Kan extension of
$G$   along $\gamma '$  (see \cite{ML}(X.3)), so $G$ is left
derivable with respect to $\gamma'$. Moreover, the canonical
morphism
$$\theta_{\gamma',G}:(Ran_{\gamma'}G)\gamma'=G \lambda\gamma' \lra
G$$ is defined by $G(\varepsilon'_{ X})$, for each object $X$ of
$\C[\mathcal S^{-1}]$.

By Lemma \ref{lema1}, $F'=Ran_\delta F$ and
$\theta_{\delta,F}=\id$. Since $Ran_{\gamma'}F'=F' \lambda$ we
have, by Lemma \ref{lema2} below, $$Ran_\gamma
F=Ran_{\gamma'}(Ran_\delta F)=F' \lambda$$ so $\mathbb L_\W
F=(Ran_\gamma F) \gamma=F' \lambda \gamma$. In addition, for each
object $X$, the canonical morphism $(\theta_{\gamma, F})_X$ is
defined by
$$(\theta_{\gamma,F})_X=(\theta_{\gamma',F'})_{\delta X}\circ
(\theta_{\delta,F})_{X}=F'(\varepsilon'_{\delta X}) .$$
\end{proof}

\begin{lemma}\label{lema2}
Let $\gamma_1:\C_1\lra \C_2$ and $\gamma_2:\C_2\lra \C_3$ be two
composable functors, and $\gamma=\gamma_2\gamma_1$. If $F:\C_1\lra
\D$ is a functor such that $Ran_{\gamma_2} (Ran_{\gamma_1}(F)) $
exists, then
\begin{enumerate}[(1)]
\item
$Ran_\gamma F$ exists, $Ran_\gamma F=Ran_{\gamma_2}
(Ran_{\gamma_1}(F))$; and
\item  $\theta_{\gamma,F}=\theta_2\gamma_1\circ \theta_1$, where
$\theta_2=\theta_{\gamma_2,Ran_{\gamma_1}(F)}$ and
$\theta_1=\theta_{\gamma_1,F}$.  \hfill{$\Box$}
\end{enumerate}
\end{lemma}

\begin{example} The previous theorem is an extension to
a non-necessarily  additive setting of the standard derivability
criterion for additive functors (see \cite{GM}, III.6, th. 8).  In
fact, let $\A$ and $\B$ be abelian categories. Suppose that $\A$
has enough projective objects, hence, by Example \ref{EjemploCE},
$(\mathbf C_+(\A),\mathcal S,\W)$ is a left Cartan-Eilenberg
category. Let $F:\mathbf C_{+}(\A)\lra \mathbf K_{+}(\B)$ be a
functor induced by  an additive functor $\A\lra \B$. Then, since
$F$ is additive, it sends homotopy equivalences to isomorphisms,
hence, by Theorem \ref{derivarconS},  $F$  is  left derivable and
$\mathbb L_\W F=F'\circ \lambda\circ \gamma$.
\end{example}

Next we study the Cartan-Eilenberg structure on the category
$\Cat_{\mathcal S}(\C,\D)$.
\begin{theorem}\label{derivarconS}
Let $(\C,\mathcal S,\W)$ be a left Cartan-Eilenberg category and
$\D$ any category. Consider the category with weak equivalences
$(\Cat_{\mathcal S}(\C,\D),\widetilde{\W})$, where
$\widetilde{\W}$ is the class of morphisms of functors
$\phi:F\Rightarrow G:\C\lra \D$ such that $\phi_M$ is an
isomorphism for all cofibrant objects $M$ of $\C$. The functor
$$ \mathbb L_\W:\Cat_{\mathcal S}(\C,\D)\lra \Cat_{\mathcal
S}(\C,\D), \quad  \mathbb L_\W F:=F' \lambda \gamma,$$ together
with the natural transformation $\theta:\mathbb L_\W F\Rightarrow
F $ defined by $(\theta_F)_X=F'(\varepsilon'_{\delta(X)}),$ for
each object $X$ of $\C$, satisfy

\begin{enumerate}
\item $(\mathbb L_\W,\theta)$ is a left resolvent functor on
$(\Cat_{\mathcal S}(\C,\D),\widetilde{\W})$;

\item  $(\Cat_{\mathcal
S}(\C,\D),\widetilde{\W})$ is a left Cartan-Eilenberg category;
and

\item  $\Cat_\W(\C,\D)$ is the subcategory of its cofibrant objects.

\end{enumerate}
\end{theorem}
\begin{proof} Since $\mathcal S\subset \W$, the category
$\Cat_\mathcal S(\C,\D)$ contains
$\Cat_\mathcal W(\C,\D)$ as a full subcategory.  On the other
hand, by Theorem \ref{derivarS}, $\Cat_\mathcal S(\C,\D)$  is a
full subcategory of $\Cat'((\C,\W),\D)$. Therefore, by Theorem
\ref{derivargenerico}, $(\mathbb L,\theta)$ induces a coidempotent
functor on $\Cat_\mathcal S(\C,\D)$, whose essential  image is
$\Cat_\mathcal W(\C,\D)$. In addition, by Theorem \ref{MA},
$(\Cat_{\mathcal S}(\C,\D),\widetilde{\W})$ is a left
Cartan-Eilenberg category whose cofibrant objects are functors in
$\Cat_\mathcal W(\C,\D)$, and $(\mathbb L,\theta)$ is a resolvent
functor, where $(\theta_F)_X=F'(\varepsilon'_{\delta(X)})$, by
Theorem \ref{derivargenerico}.

Next, by Theorem \ref{derivarS}, $\mathbb L F=F'\lambda\gamma$
and, by Theorem \ref{MA}, the class of weak equivalence is the
class of morphisms $\phi:F\Rightarrow G$ such that $\mathbb
L(\phi)$ is an isomorphism, that is to say,
$\phi_{\lambda(\gamma(X))}$ is an isomorphism, for each $X$. Since
the objects $\lambda(\gamma(X))$ are  the cofibrant objects up to
strong equivalences, a morphism $\phi$ is a weak equivalence if
and only if $\phi_M$ is an isomorphism for each cofibrant object
$M$, that is to say, $\widetilde{\W}$ is the class of weak
equivalences.
\end{proof}

\subsection{Models of  functors}\label{section3.3}

When the target category $\D$ of  functors $F:\C\lra \D$ is
endowed with a class of weak equivalences $\mathcal E$, the
previous results can be applied to the functor $\gamma_\D F:\C\lra
\D[\E^{-1}]$ to obtain a model of this functor. However, in some
situations, it is desirable to have cofibrant models for the
functor $F$ itself.
  We prove that this is possible if $\C$ is a left
Cartan-Eilenberg category with a left resolvent functor and $F$
sends strong equivalences to weak equivalences.

\begin{nada}\label{notacion3.3}
So let $(\C,\mathcal S,\W)$ be a Cartan-Eilenberg category with a
left resolvent functor $(R,\varepsilon)$ and $\D$  a category with
a saturated class of weak equivalences $\mathcal E$. Denote by $
\Cat_{\mathcal S,\E}(\C,\D) $ the full subcategory of
$\Cat(\C,\D)$ whose objects are the functors  which send $\mathcal
S$ to $\mathcal E$.
\end{nada}
\begin{definition}
Let  $F,G$ be objects of $\Cat_{\mathcal S,\mathcal E}(\C,\D)$ and
$\phi:F\Rightarrow G$ a morphism.
\begin{itemize}
\item[(i)] $\phi $ is called a \emph{weak equivalence} if $\phi_M$
is in $\mathcal E$, for all $M\in Ob\ \C_{cof}$. \item[(ii)] $\phi
$ is called a \emph{strong equivalence} if $\phi_X$ is in
$\mathcal E$, for all $X\in Ob\ \C$.
\end{itemize}
\end{definition}

We denote by $\widetilde{\W}$  and $\widetilde{\mathcal S}$ the
classes of  weak and strong equivalences of $\Cat_{\mathcal
S,\mathcal E}(\C,\D)$, respectively.

If $F(\mathcal S)\subset \mathcal E$, then $R^*(F)(\mathcal
S)=F(R(\mathcal S))\subset F(\overline{\mathcal S}) \subset
\mathcal E$, thus the resolvent functor $R$ induces the functor
$$R^*:\Cat_{\mathcal S,\mathcal E}(\C,\D)\lra \Cat_{\mathcal
S,\mathcal E}(\C,\D)$$ given by ${R^*}(F):=FR$, and the counit
$\varepsilon:F\Rightarrow \id$ induces a counit
${\varepsilon^*}:{R^*}\Rightarrow \id$ by
$$\varepsilon^*_F:=F\varepsilon:F R\longrightarrow F.$$

\begin{theorem}\label{Resolventeparafuntores} Let $(\C,\mathcal S,\W)$
be a  category with a
left resolvent functor  $(R,\varepsilon)$,  and $\D$  a category
with a saturated class of weak equivalences $\mathcal E$. With the
previous notation we have
\begin{enumerate}[(1)]
\item  $({R^*},{\varepsilon^*})$ is a
left resolvent functor for $\left(\Cat_{\mathcal S,\mathcal
E}(\C,\D),\widetilde{\mathcal S}, \widetilde{\W}\right)$;

\item  $\left(\Cat_{\mathcal S,\mathcal E}(\C,\D),\widetilde{\mathcal S},
\widetilde{\W}\right)$ is a left Cartan-Eilenberg category; and

\item a functor  $F\in Ob\, \Cat_{\mathcal S,\mathcal E}(\C,\D)$ is cofibrant
if and only if  $F(\W)\subset\mathcal E$.
\end{enumerate}
\end{theorem}

\begin{proof}
We first observe that, by (2) of Lemma \ref{LFR}, for each object
$X$ of $\C$, $\varepsilon_{RX}$ and $R(\varepsilon_X)$ are in
$\overline{\mathcal S}$, therefore, for each functor $F$ in
$\Cat_{\mathcal S,\mathcal E}(\C,\D)$, the morphisms
$F(\varepsilon_{RX})$ and $F(R(\varepsilon_X))$ are in $\mathcal
E$, hence $R^*\varepsilon_F$ and $\varepsilon_{R^*(F)}$ are in
$\widetilde{\mathcal S}$.

Moreover, by (3) of Proposition \ref{PreMA}, it is easy to check
that $\widetilde{\W}=\left({R^*}\right)^{-1}(\widetilde{\mathcal
S})$. In particular ${R^*}(\widetilde{\mathcal S})\subset
\widetilde{\mathcal S}$. Hence we can  apply Theorem \ref{MA} to
obtain (1) and (2).

By part (1) and Proposition \ref{PreMA}, $F$ is cofibrant if and
only if $\varepsilon^*_F:R^*(F)\lra F$ is a strong equivalence,
that is to say, $F(\varepsilon_X):F(RX)\lra F(X)\in \mathcal E$,
for each $X$.

Hence, if $F(\W)\subset \mathcal E$, since $\varepsilon_X\in \W$,
we obtain $F(\varepsilon_X)\in \mathcal E$, that is to say
$\varepsilon^*_F:R^*F\lra F$ is a strong equivalence.

To prove the converse, observe that if $F$ is a functor such that
$F(\mathcal S)\subset \mathcal E$, then  we have also
$F(\overline{\mathcal S})\subset \mathcal E$ since $\mathcal E$ is
saturated. By Lemma \ref{LFR}, for each $w\in \W$, we have
$R(w)\in \overline{\mathcal S}$, so $F(R(w))\in \mathcal E$. Hence
$F(R(\W))\subset \mathcal E$.

Now, suppose that $F$ is cofibrant, and let $w:X\lra Y\in \W$. We
have a commutative diagram
$$
\xymatrix{
FRX\ar[dd]_{FRw}\ar[rr]^{F\varepsilon_X}&&F(X)\ar[dd]^{Fw}
\\ \\
 FRY\ar[rr]^{F\varepsilon_Y}&&F(Y)
}.
$$ Since $F\varepsilon_X$, $F\varepsilon_Y$ and $FRw$ are in $\mathcal E$,
we obtain $Fw\in \mathcal E$, since $\mathcal E$ is saturated,
that is to say $F(\W)\subset \mathcal E$.
\end{proof}

Finally, by Theorems  \ref{Resolventeparafuntores} and
\ref{derivarS}, we obtain:

\begin{corollary}\label{DerivedFunctor2}
 With the previous notation, for each $F\in
\Cat_{\mathcal S,\mathcal E}(\C,\D)$, $F\varepsilon: F R \lra F$
is a cofibrant left model of $F$,  the left derived functor
$\mathbb L_\W (\gamma_\E F)$ of $\gamma_\E F $  is $\gamma _\E F
R$, and the total left derived functor $\mathbb L F$  of $F$ (see
\cite{Q}, Definition 2, \S I.4) is  the functor induced by
$\mathbb L_\W( \gamma_\E F)$, so we have a commutative diagram
$$
\xymatrix{ \C\ar[ddrr]_{\mathbb L_{\W}(\gamma_\E
F)}\ar[dd]_{\gamma_{\W} }\ar[rr]^{F\circ R}
&&\D\ar[dd]^{\gamma_\E}\\ \\
\C[\W^{-1}]\ar[rr]^{\mathbb L F}&&\D[\E^{-1}]\,. }
$$
\end{corollary}

\begin{example}
 Let  $\mathbf C_+(A)$ be the Cartan-Eilenberg
category  of bounded below chain complexes of $A$-modules, where
$A$ is a commutative ring,  and $\varepsilon:R\Rightarrow \id$ the
resolvent functor defined by the free functorial resolution (see
Example \ref{ComplejosSonCE}). Let $\B$ be an abelian category and
 $F:\mathbf C_+(A)\lra \mathbf C_+(\B)$ a functor induced by an additive functor $A-\text{mod}
\lra \B$. Then $F$  sends homotopy equivalences to
quasi-isomorphisms, therefore $F\varepsilon:F R\Rightarrow F$ is a
cofibrant left model of $F$ in $\Cat_{\mathcal S,\E}(\mathbf
C_+(A), \mathbf C_+(B))$, where $\mathcal S$ are the homotopy
equivalences and $\E$ the quasi-isomorphisms.
\end{example}

\section{Quillen model categories and Sullivan minimal models}
\label{Quillenmodels}

In this section we describe how Cartan-Eilenberg categories relate
to some other axiomatisations for homotopy theory.

\subsection{Quillen model categories} Let $\C$ be a Quillen model
category, that is, a category equipped with three classes of
morphisms: weak equivalences $\W$, cofibrations $cofib$, and
fibrations $fib$, satisfying Quillen's axioms for a  model
category (\cite{Q}).

In a Quillen model category there are the notions of cofibrant,
fibrant and cylinder objects. To distinguish between these objects
and the cofibrant/fibrant/cylinder objects as introduced in this
paper, the former ones will be called Quillen
cofibrant/fibrant/cylinder objects. Denote by $\C_f$ and $\C_{cf}$
the full subcategories of Quillen fibrant and cofibrant-fibrant
objects of $\C$, respectively.

In a Quillen model category there are the notions of left and
right homotopy. For instance, if $f,g:X\lra Y$ are two morphisms,
a left homotopy from $f$ to $g$ is a morphism $h:X'\lra Y$, where
$X' $ is a Quillen cylinder object for $X$ (that is,
$\partial_0\vee\partial_1:X\vee X\lra X'$ is a cofibration,
$p:X'\lra X$ is a weak equivalence, and $p\partial_0=Id=p\partial
_1$, see  Definition I.4 of \cite{Q}), such that $h\partial _0=f$
and $h\partial _1=g$. Let $\sim_l$ be the equivalence relation
transitively generated by the left homotopy, and let $\mathcal
S_l$ be the class of homotopy equivalences coming from $\sim_l$.
We denote by $\pi^l(X,Y)$ the set of equivalence classes of
morphisms from $X$ to $Y$ with respect to $\sim_l$. By the dual of
(\cite{Q}, Lemma I.6), $\sim_l$ is a congruence in $\C_f$.

\begin{lemma}
The equivalence relation $\sim_l$ is  compatible with $\mathcal
S_l$ in $\C_f$.
\end{lemma}
\begin{proof}
Let $f,g:X\lra Y$ be two  morphisms such that $f\sim_lg$, where
$X,Y$ are fibrant objects. We can assume that there exists a left
homotopy $h':X'\lra Y$ from $f$ to $g$, where $X'$ is a cylinder
object for $X$. We can choose a cylinder object such that
$p':X'\lra X$ is a trivial fibration. In fact, let
$$\xymatrix{X'\ar[r]^j& {X\times I}\ar[r]^p& X}$$ be a factorisation
of $p'$ in a trivial fibration $p$ and a cofibration $j$, which is
also trivial since $p'$ is too. Since $Y$ is a fibrant object, and
$j$ is a trivial cofibration, there exists a morphism ${h}$
filling the following solid-arrow commutative diagram.
$$\xymatrix{
X'\ar[d]^j\ar[r]^{h'}&Y\ar[d]\\
{X\times I}\ar@{.>}[ur]^{{h}}\ar[r]& {*} }$$ Therefore ${h}$ is a
left homotopy from $f$ to $g$.

Next  the trivial fibration $p:{X\times I}\lra X$ is a left
homotopy equivalence. This is a consequence of the following
general fact in a Quillen model category: If a cofibration
$i:X\lra Y$ has a retraction $p:Y\lra X$ which is a trivial
fibration, then $i$ (and $p$)  is a
 left homotopy
equivalence. (\emph{Proof}: Let $\delta_0\vee \delta_1:Y\vee Y
\lra Y\times I\lra Y$ be a Quillen cylinder object for $Y$.
Consider the diagram
$$
\xymatrix{Y\vee Y\ar[d]_{\delta_0 \vee \delta_1}
\ar[r]^{ ip\vee 1_Y}&Y\ar[d]^{p}\\
Y\times I\ar@{.>}[ru]^H\ar[r]^{pq}&X,}$$ where the left vertical
arrow is a cofibration and the right one is a trivial fibration.
Then,  the lifting $H$ is a left homotopy between $ip$ and $1_Y$.)
Going back to the proof of the lemma, since $p\in \mathcal S_l$,
we have, in $\C_f[\mathcal S^{-1}]$,
$f={h}\partial_0={h}p^{-1}p\partial_0={h}p^{-1}=h\partial_1=g$, as
asserted.
\end{proof}

By the previous lemma, the class $\mathcal S_l$ is compatible with
$\sim_l$ and, by Proposition \ref{congruencia}, there is an
isomorphism of categories $\pi^{l}\C_{f}\cong \C_f[\mathcal
S_l^{-1}]$. Therefore, the relative localisation $\C_{cf}[\mathcal
S_l^{-1},\C_f]$ is isomorphic to the homotopy category
$\pi^{l}\C_{cf}$. We observe that the left homotopy relation is,
itself, an equivalence relation when restricted to the subcategory
$\C_{cf}$, by Lemma 4 of \cite{Q}.

Let $\W$ be the class of weak equivalences of $\C_f$. Since
$p:Cyl(X)\lra X$ is a weak equivalence,  we have that $\mathcal
S_l\subset \overline{\W}$, so $(\C_f, \mathcal S_l, {\W})$ is a
category with strong and weak equivalences.

\begin{theorem}\label{exampleQuillen}
Let $\C$ be a Quillen model category. Then $(\C_{f},
\mathcal{S}_l,{\W})$ is a left Cartan-Eilenberg category and
$\C_{cf}$ is a subcategory of cofibrant left models  of $\C_{f}$.
\end{theorem}
\begin{proof}
We prove that the class $\C_{cf}$ satisfies the hypothesis of
Theorem \ref{F}. Let $M$ be a Quillen fibrant-cofibrant object,
and let $w:Y\lra X$ be a weak equivalence. Let us see that the map
$$w_*:\C[\mathcal S_l^{-1}](M,Y)=\pi^l(M,Y)\lra
\C[\mathcal S_l^{-1}](M,X)=\pi^l(M,X)$$ is  bijective. By the
axiom M2 of \cite{Q}, there exists a factorisation $w=\beta\circ
\alpha$, where $\alpha:Y\lra Z$ is a trivial cofibration and
$\beta:Z\lra X$ is a trivial fibration. Since $w_*=\beta_*\circ
\alpha_*$ it is enough to prove that the maps
$$\alpha_*:\pi^l(M,Y)\lra \pi^l(M,Z)$$ and
$$\beta_*:\pi^l(M,Z)\lra \pi^l(M,X)$$
are bijective. By Lemma 7 of \cite{Q}, $\beta_*$ is bijective.

To prove that $\alpha_*$ is also bijective, we apply the dual of
Lemma 7 of \cite{Q}, for which we denote by $\pi^r$ the right
avatar of $\pi^l$. Indeed, the map $\alpha^*:\pi^r(Y,A)\lra
\pi^r(Z,A),$ is bijective for each Quillen-fibrant object $A$, by
the dual of  Lemma 7 of \cite{Q}, since $\alpha$ is a trivial
cofibration. Therefore $\alpha$ is an isomorphism in
$\pi^r\C_{f}$, and as a consequence  the map
$\alpha_*:\pi^r(M,Y)\lra \pi^r(M,X)$ is bijective. On the other
hand,  $M$ being Quillen-cofibrant, for each Quillen-fibrant
object $X$, the left and right homotopy relations coincide in
$\C(M,X)$, hence $\alpha_*:\pi^l(M,Y)\lra \pi^l(M,Z)$ is
bijective.

Finally, by Quillen axiom M2, for each Quillen-fibrant object $X$
there exist a trivial fibration $M\lra X$, where $M$ is
Quillen-cofibrant, and moreover $M$ is Quillen fibrant, by M3.
\end{proof}

\begin{remark}
Observe that in a Quillen model category $\C$ the definition of
Quillen cofibrant objects is not homotopy invariant, while the
subcategory of  cofibrant objects of $\C_f$ is stable by homotopy
equivalences. In fact, the  cofibrant objects are those homotopy
equivalent to Quillen cofibrant objects.

For instance, let $\A$ be an abelian category with enough
projectives and $\PChains{\A}$ the category of {bounded below}
chain complexes. It is well known (see \cite{Q}, Chapter I)  that
taking quasi-isomorphisms as weak equivalences, epimorphisms as
fibrations, and monomorphisms whose cokernel is a degree-wise
projective complex as cofibrations, $\PChains{\A}$ is a Quillen
model category with all objects fibrant. A contractible complex is
cofibrant, but it is not Quillen cofibrant unless it is projective
(see also \cite{C}).
\end{remark}

\subsection{Sullivan minimal models} In some Cartan-Eilenberg
categories there is a distinguished subcategory $\M$ of $\C_{cof}$
which serves as a subcategory of cofibrant left models. A typical
situation is that of Sullivan minimal models (\cite{S}). Let us
give an abstract version.

\begin{definition}\label{defminimal}
Let $(\C,{\mathcal S}, \W)$ be a category with strong and weak
equivalences. We say that a cofibrant object $M$ of $\C$ is
\emph{minimal} if
$$
End_\C(M)\cap \W = Aut_\C(M),
$$
that is, if any weak equivalence $w:M\lra M$ of $\C$ is an
isomorphism.
\end{definition}

We denote by $\C_{min}$ the full subcategory of $\C$ whose objects
are minimal in $(\C, \mathcal{S}, \W)$.

\begin{definition}
We say that $(\C, \mathcal{S}, \W)$ is a \emph{left Sullivan
category} if there are enough minimal left models.
\end{definition}

\begin{remark}\label{minimalnohomotopia}
Observe that by the uniqueness property of the extension in
Definition \ref{defcofobj}, any cofibrant object of $\C$ is
minimal in the localised category $(\loc{C}{S}, \delta(\W))$.
\end{remark}

\begin{remark}\label{minimal no localiza}
As a consequence of the definition, a left Sullivan category is a
special kind of a left Cartan-Eilenberg category, one for which
the canonical functor
$$
\C_{min}[\mathcal S^{-1}, \C]\lra \C[\W^{-1}]
$$
 is an
equivalence of categories. Observe that by definition, if $X$ is a
minimal object and $s:X\lra X$ is in $\mathcal{S}$, then $s$ is an
isomorphism, hence $\C_{min}[\mathcal S^{-1}]=\C_{min}$, so that
in this case the inclusion functor $\C_{min}[\mathcal S^{-1}]\lra
\C_{min}[\mathcal S^{-1},\C]$ is not, generally speaking, an
equivalence of categories.

\end{remark}

\begin{nada} An example of a Sullivan category is provided by
the original Sullivan's minimal cdg algebras. Let $\mathbf k$ be a
field of characteristic zero, and $\Adgc_1$ the category of simply
connected commutative differential graded $\mathbf k$-algebras
($1$-connected $\mathbf k$-cdg algebra, for short).

A path object for a $\mathbf k$-cdg algebra $B$ is the tensor
product $\Path(B):=B\otimes \mathbf k[t,dt]$, together with the
morphisms $\delta_0,\delta_1:\Path(B)\lra B$, and $p:B\lra
\Path(B)$ defined by $\delta_i(a(t))=a(i)$ for $i=0,1$, and
$p(a)=a\otimes 1$.

Let $f_0,f_1:A\lra B$ be two morphisms of $\mathbf k$-cdg
algebras. A right homotopy from $f_0$ to $f_1$ is a morphism of
$\mathbf k$-cdg algebras, $H:A\lra \Path(B)$ such that $\delta_i
H=f_i, i =0,1$ (see \cite{S} or \cite{GM}, (10.1)).

Let $\sim$ be the equivalence relation transitively generated by
the right homotopy. It follows from the functoriality of the path
object that $\sim$ is a congruence. Let $\mathcal S$ be the class
of homotopy equivalences with respect to $\sim$.
\end{nada}

\begin{lemma}\label{compatibleS}
The equivalence relation $\sim$ is  compatible with $\mathcal S$.
\end{lemma}

\begin{proof} Because of Example \ref{cilindro}, it is enough to see
that $p:B\lra \Path(B)$ is in $\mathcal S$ and this follows from
the fact that $\delta_0 p=\id_B$ and $H:\Path(B)=B\otimes \mathbf
k[t,dt]\lra \Path(\Path(B))=(B\otimes \mathbf k[t,dt])\otimes
\mathbf k[u,du]$ defined by $H(a(t))=a(tu)$ is a right homotopy
from $p\delta_0 $ to $Id_{\Path(B)}$.
\end{proof}

So, by  Proposition \ref{congruencia}, there is an isomorphism of
categories $$\Adgc_1/\!\!\sim \,\, \cong \Adgc_1[\mathcal
S^{-1}].$$

Let $\W$ be the class of quasi-isomorphisms of $\Adgc_1$; that is,
those morphisms inducing isomorphisms in cohomology. Since
$p:B\lra \Path(B)$ is a quasi-isomorphism, we have that $\mathcal
S\subset \overline{\W}$. So $(\Adgc_1, \mathcal S, {\W})$ is a
category with strong and weak equivalences.

Recall that a $\mathbf k$-cdg algebra $A$ is a $1$-connected
Sullivan minimal $\mathbf k$-cdg algebra if it is a free graded
commutative $\mathbf k$-algebra $A=\Lambda(V)$ such that
$A^0=\mathbf k$, $A^1=0$, and  $dA^+\subset A^+\cdot A^+$, where
$A^+=\oplus_{i>0} A^{i}$ (\cite{S}, see also \cite{GM}, p. 112).
Let $\M_S$ be the full subcateogory of $1$-connected Sullivan
minimal $\mk$-cdg algebras. We can sum up Sullivan's results on
minimal models in the following theorem.

\begin{theorem} \label{adgc}
$(\Adgc_1,\mathcal S,\W)$  is a left Sullivan category and $\M_S$
is the subcategory of minimal objects of $\Adgc_1$.
\end{theorem}

\begin{proof}
First of all, let us check the hypotheses of Theorem \ref{F} for
the class $\M_S$ of Sullivan minimal $1$-connected algebras. Let
$M$ be a $1$-connected Sullivan minimal $\mathbf k$-cdg algebra.
If $A\lra B$ is a quasi-isomorphism, the induced map $[M,A]\lra
[M,B]$ between the sets of homotopy classes of morphisms is
bijective, by \cite{GM} Theorem 10.8. So $M$ is a cofibrant
object. In addition, by \cite{GM} Theorem 9.5, any $1$-connected
$\mathbf k$-cdg algebra has a Sullivan minimal model, so, by
Theorem \ref{F}, $\M$ is a subcategory of left cofibrant models of
$\Adgc_1$.

By \cite{GM} Lemma 10.10, any quasi-isomorphism $M\lra M$ of a
Sullivan minimal algebra is an isomorphism, so $M$ is a minimal
object in $(\Adgc_1,\mathcal S,\W)$, therefore $(\Adgc_1,\mathcal
S,\W)$  is a left Sullivan category.

Reciprocally, every minimal object of $\Adgc_1$ is isomorphic to a
Sullivan minimal $1$-connected algebra. Let $M$ be a minimal
object of $\Adgc_1$. Because of \cite{GM} Theorem 9.5, there is a
Sullivan minimal model $\omega : M_S \longrightarrow M \in \W$.
Since $M$ is a cofibrant object, we have a bijection $\omega_* :
[M,M_S] \longrightarrow [M,M]$. Let $\phi : M \longrightarrow M_S$
be such that $\omega\phi \sim \id_M$. Then $H(\omega\phi ) =
\id_{HM}$ and so $\omega\phi$ is an isomorphism, because $M$ is a
minimal object. Also because of the \emph{2 out of 3} property of
quasi-isomorphisms, $\phi \in \W$. So again we find $\psi : M_S
\longrightarrow M$ such that $\phi\psi \sim \id_{M_S}
\Longrightarrow \psi \sim \omega $, which also implies that
$\phi\omega \sim \id_{M_S}$. So, $\phi\omega$ is an isomorphism
too. Hence so is $\omega$.
\end{proof}

\begin{nada}
Analogously, there are enough minimal objects in the category $\Op
(\mk)_1$ of dg operads over $\mathbf k$, $P$, such that $H^*P(1) =
0$, (see \cite{MSS}). From Theorem \ref{F} again it follows that
$\Op (\mk)_1$ is a left Sullivan category.

We next consider in greater detail the case of dg modular operads
over a field of characteristic zero $\mathbf k$ (refer to
\cite{GK} and \cite{GNPR1} for the notions concerning modular
operads that will be used).

Let $\MOp (\mk)$ be the category of dg modular operads. We have an
analogous path object for modular operads: if $P$ is a dg modular
operad, its path object is the tensor product $\Path (P) = P
\otimes \mk [t, \delta t]$. Let $\sim$ be the equivalence relation
transitively generated by the right homotopy defined with this
path object. We can see, as in Lemma \ref{compatibleS}, that the
class of homotopy equivalences $\Scal$ with respect to $\sim$ is
compatible with $\sim$, so we have an isomorphism of categories
$\MOp (\mk)/\!\!\sim \,\,  \cong \MOp (\mk)[{\Scal}^{-1}]$.

Let $\W$ be the class of quasi-isomorphisms of $\MOp (\mk)$. We
see in the same way as for $\Adgc_1$ that $(\MOp (\mk), \Scal ,
\W)$ is a category with strong and weak equivalences.

In \cite{GNPR1}, Definition 8.6.1, we defined minimal modular
operads as modular operads obtained from the trivial operad $0$ by
a sequence of principal extensions. Let $\M$ be the full
subcateogory of minimal modular operads.
\end{nada}

\begin{theorem}\label{MOP}
$(\MOp (\mk), \Scal , \W)$ is a left Sullivan category and $\M$ is
the subcategory of minimal objects of $\MOp (\mk)$.
\end{theorem}

\begin{proof} Let us check the hypothesis of Theorem \ref{F}: if
$M$ is a minimal modular operad and $P \longrightarrow Q$ a
quasi-isomorphism of $\MOp (\mk)$, the induced map $[M,P]
\longrightarrow [M,Q]$ is a bijection by \cite{GNPR1}, Theorem
8.7.2. So $M$ is a cofibrant object. The existence of enough
cofibrant objects is guaranteed by Theorem 8.6.3. op.cit., and
these minimal modular operads are minimal objects because of
op.cit., Proposition 8.6.2.

We can argue as in the proof of Theorem \ref{adgc} to show that
every minimal object of $\MOp (\mk)$ is isomorphic to an object of
$\M$.
\end{proof}

\section{Cartan-Eilenberg categories of filtered objects}

In this section we prove that some categories of filtered
complexes and of filtered graded differential commutative algebras
are Cartan-Eilenberg categories.

\subsection{Filtered complexes of an abelian category}
Let $\A$ be an abelian category. By a filtered complex of $\A$ we
understand a pair $(X,W)$ where $X$ is a chain complex of $\A$ and
$W$ is an increasing filtration of $X$ by subcomplexes  $W_pX$. We
denote by $Gr^W_pX$ the complex $W_pX/W_{p-1}X$.

\begin{nada}
We denote by $\FPChains{\A}$ the  category of filtered complexes
$(X,W)$ such that
\begin{enumerate} [(i)]
\item the complex $X$ is
bounded below, that is, $X_p=0$ if $p\ll 0$;
\item the filtration $W$ is bounded below and
biregular, that is $W_pX=0$ if $p\ll 0$
 and $W$
is finite on each $X_n$.
\end{enumerate}
\end{nada}

\begin{nada}
We recall that, given a chain complex $X$ and an integer $n$,
$X[n]$ denotes the  chain complex which in degree $i$ is equal to
$X_{i-n}$ with differential $(-1)^n\partial _{i-n}$. If, in
addition, $W$ is a filtration on $X$, then $X[n]$ has an induced
filtration defined by $W_p(X[n])=(W_pX)[n]$.

Two morphisms $f,g:(X,W)\lra (Y,W)$ between filtered chain
complexes are  \emph{filtered homotopic} if there is a
\emph{filtered homotopy} from $f$ to $g$, that is,  a homotopy
$h:X\lra Y$ of chain maps from $f$ to $g$ which is a filtered
homogeneous morphism. The filtered homotopy relation is an
equivalence relation  $\sim$ which is compatible with composition,
hence is a congruence on $\FPChains{\A}$.
 We
denote by $\mathcal{S}$ the class of filtered homotopy
equivalences, and by $[(X,W),(Y,W)]$ the set of homotopy classes
of  morphisms $(X,W)\lra (Y,W)$. The quotient category $\mathbf
{KF}_+(\A):= \FPChains{\A}/\sim$ is called the \emph{filtered
homotopy category} of $\A$.

The filtered homotopy relation is compatible with the class
$\mathcal S$, since it is easy to prove that this relation can be
expressed by a cylinder object. Recall that, given a filtered
chain complex $X$, the cylinder of $X$ is the filtered complex
$Cyl(X) =X \oplus X[1] \oplus X $ with differential
$$
D = \left(
\begin{array}{ccc}
\partial  & - \id  &  0   \\
 0&  -\partial  &  0   \\
 0& \id  &  \partial
\end{array} \right),
$$
together with  the two injections, $i_0,i_1:X\lra \Cyl(X)$,  and
the filtered homotopy equivalence $p:\Cyl(X)\lra X$, defined by
$\begin{pmatrix}\id&0&\id\end{pmatrix}$, (cf. \cite{GMa}). Hence,
by Proposition \ref{congruencia}, the categories $\mathbf
{KF}_+(\A)$ and $\FPChains{\A}[\mathcal S^{-1}]$ are canonically
isomorphic.

\end{nada}

A filtered morphism $f$ is called a \emph{filtered
quasi-isomorphism}  if $W_p(f)$ is a quasi-isomorphism for each
$p$ (equivalently, since filtrations are bounded below and
biregular, if $Gr^W_p(f)$ is a quasi-isomorphism for each $p$).
Denote by $\W$ the class of filtered quasi-isomorphisms in
$\FPChains{\A}$. The localised category
$\mathbf{FC}_+(\A)[\W^{-1}]$ is the filtered derived category of
$\A$, $\mathbf {DF}_+(\A)$ (see \cite{I}).

It is clear that $\mathcal{S}\subset \W$, so
$(\FPChains{\A},\mathcal{S},\W)$ is a category with strong and
weak equivalences.

\begin{theorem}\label{CEfiltrados}
Let $\A$ be an abelian category with enough projective objects,
and let $\mathcal P$ the full subcategory of  filtered complexes
$P$ such that, for all $p$, $Gr^W_pP$ is projective in each
degree. Then,
\begin{enumerate}[(1)]
\item every object in $\mathcal P$ is cofibrant;

\item the category $(\FPChains{\A},\mathcal{S},\W)$ is a left
Cartan-Eilenberg category; and

\item the functor $\mathcal P/\sim\, \lra \mathbf{DF}_+(\A)$ is an
equivalence of categories.
\end{enumerate}
\end{theorem}

The theorem follows from Theorem \ref{F}, as soon as we check its
hypotheses in Propositions \ref{i-filtrado},
\ref{ii-filtrados} and \ref{fil2} below.

\begin{proposition}\label{i-filtrado}
For any filtered quasi-isomorphism $w:Y\lra X$, and any morphism
$f:P\lra X$ in $\FPChains{\A}$, where $P\in Ob\, \mathcal P$,
there exists a morphism $g:P\lra Y$ in  $\FPChains{\A}$ such that
$wg$ is filtered homotopic to $f$.
\end{proposition}
\begin{proof}Let
$$\xymatrix{
&Y\ar[d]^w\\P\ar@{.>}[ur]^g\ar[r]^f&X }$$

 be a solid diagram in  $\FPChains{\A}$, where  $w$ is a filtered quasi-isomorphism. To prove
the existence of a  lifting of $f$, up to homotopy, we define,
inductively on $p$, a chain map $g_p:W_pP\longrightarrow W_pY$ and
a homotopy $h_p:W_pP\longrightarrow W_pX$ such that $h_p:
wg_p\stackrel{\sim}\ra f$, as follows. For $p\ll 0$ we have
$W_pP=0$, so we take $g_p=0$ and $h_p=0$. Assume now that
$g_{p-1}$ and $h_{p-1}$ have been defined, and consider the
diagram
$$
\begin{CD}
W_{p-1}P @>g_{p-1}>>  W_pY\\@Vj_pVV@VVW_pwV\\
W_{p}P   @>f|W_pP>>      W_pX\ .
\end{CD}
$$
The cokernel of $j_p$ is projective in each degree and bounded
below, and $W_{p}w$ is a quasi-isomorphism, hence, by the lifting
up to homotopy property (which we will recall in the next Lemma
\ref{fillema}), there are a chain map $g_p:W_pP\lra W_pY$ and a
homotopy $h_p:W_pP\lra W_pX$, which are extensions of the previous
data $g_{p-1}$ and $h_{p-1}$. As the filtration $W$ is biregular,
$g_p$ and $h_p$ define a filtered morphism $g:P\longrightarrow Y$
and a filtered homotopy $h:P\longrightarrow X$, respectively, such
that $h: wg\stackrel{\sim}\ra  f$.
\end{proof}

\begin{lemma}\label{fillema} { \rm(Lifting up to homotopy
property.)} Let
$$\xymatrix{
Q\ar[dd]_j\ar[rr]^\phi&&Y\ar[dd]^w\\
\\R\ar@{.>}[uurr]^G\ar[rr]^F&&X }$$ be a  diagram of chain
complexes of $\A$ such that $w$ is a quasi-isomorphism, $j$ is a
monomorphism such that $\text{\rm Coker } j$ is bounded below and
projective in each degree, and there is a homotopy   $\lambda:
w\phi\stackrel{\sim}\ra Fj $. Then, there is a chain map $G:R\lra
Y$ such that $G j=\phi$, and a homotopy $H:w G\stackrel{\sim}\ra F
$, such that $Hj=\lambda$.
\end{lemma}

\begin{proof} Firstly,
   we recall some standard
facts about the mapping cone of a morphism and the complex of
homomorphism between two complexes (see \cite{GMa}).

 Given a morphism of complexes $f:B\lra A$, we denote by
 $C(f)$ the \emph{mapping cone} of $f$,   defined as  $C(f)=
B[1]\oplus A$ with  differential
$$
\partial^{C(f)}=\begin{pmatrix}-\partial^B&0\\
f&\partial^A
\end{pmatrix}.
$$

Given two complexes $A$ and $B$, the \emph{complex of homogeneous
homomorphisms},  $\underline{\Hom}_*(A,B)$, is defined by
$\underline{\Hom}_n(A,B)=\prod_i\Hom(A_i,B_{i+n})$, with
differential $D$ given by $Df=\partial f-(-1)^{|f|}f\partial$,
where $|f|$ is the degree of the homogeneous morphism $f$. In
particular, if $f$ is of degree $0$, then $Df=0$ if and only if
$f$ is a morphism of complexes. In addition, given two morphisms
$f,g:A\lra B$, a homogeneous morphism  $h\in
\underline{\Hom}_1(A,B)$ is a homotopy $h:f\stackrel{\sim}\ra g$
if and only if $Dh=g-f$. As a consequence, $[A,B]\cong
H_0\underline{\Hom}_*(A,B)$.

Let $P:=\text{Coker} j$. As  $P$ is projective in each degree, the
exact sequence of complexes
$$
\xymatrix{0\ar[r]&Q\ar[r]^{j}&R\ar[r]&P\ar[r]&0}
$$
 splits degree-wise, therefore we have a morphism of  exact
sequences
$$
\xymatrix{ 0\ar[r]&\underline{\Hom}_*(P,Y)\ar[dd]^{w_*^P}\ar[r]
&\underline{\Hom}_*(R,Y)\ar[r]^{j^*}\ar[dd]^{w_*^{R}}
&\underline{\Hom}_*(Q,Y)\ar[dd]^{w_*^Q}\ar[r]&0\\
\\
0\ar[r]&\underline{\Hom}_*(P,X)\ar[r]
&\underline{\Hom}_*(R,X)\ar[r]^{j^*}&\underline{\Hom}_*(Q,X)\ar[r]&0.}
$$

Since  $P$ is bounded below and projective component-wise and $w$
is a quasi-isomorphism, it is a well known  fact that $w_*^P$ is a
quasi-isomorphism,  hence its cone $C(w_*^P)$ is an acyclic
complex and, therefore, the  epimorphism
$$
j^*:C(w_*^R)\lra C(w_*^Q)
$$
is a quasi-isomorphism.

We observe that $(\phi,\lambda)\in \underline{Hom}_0(Q,Y)\oplus
\underline{Hom}_1(Q,X)=C(w_*^Q)_1$ and $(0,F)\in C(w_*^R)_0$
satisfy
$$\partial (\phi,\lambda)=(-D\phi,w^Q_*\phi+D\lambda )=(0,Fj)=j^*(0,F),\quad
\partial(0,F)=0.$$

Then we can use  the following elementary fact about complexes of
abelian groups. If $f:B\lra A$ is an epimorphism of complexes of
abelian groups which  is a quasi-isomorphism, given $a\in A_{1}$
and $b\in B_{0}$ such that $\partial a=f(b)$ and $\partial b=0$,
there exists $c\in B_{1}$ such that $f(c)=a$ and $\partial c=b$.
Indeed, since $f$ is surjective, there exists $c_0\in B_{1}$ such
that $f(c_0)=a$. Then $f(b-\partial c_0)=fb-f\partial c_0=\partial
a-\partial fc_0=0$, and $\partial (b-\partial c)=0$, hence
$b-\partial c_0\in (\text{Ker }f)_{0}$ is a cycle. Since
$\text{Ker } f$ is an acyclic complex, there exists $h\in
(\text{Ker } f)_{1}$ such that $\partial h=b-\partial c_0$.
Therefore $c:=c_0+h\in B_{1}$ satisfies $f(c)=f(c_0)+f(h)=a$ and
$\partial c=\partial c_0+\partial h=b$.

Coming back to the proof of the lemma, by the previous assertion,
there exists $(G,H)\in C(w^R_*)_1$ such that
$$j^*(G,H)=(\phi,\lambda), \quad \partial (G,H)=(0,F).$$ This means
that  $j^*G=\phi$, $j^*H=\lambda$, $DG=0$ and $w_*G+DH=F$, that
is, $G:R\lra Y$ is a chain map such that $Gj=\phi$, and
$H:F\stackrel{\sim}\ra wG$ is a homotopy such that $Hj=\lambda$,
hence $(G,H)$ verifies the  statement.
\end{proof}

\begin{proposition}\label{ii-filtrados}
For any filtered quasi-isomorphism $w:Y\longrightarrow X$  and any
$P\in Ob\,\mathcal P$, the map $w_*:[(P,W),(Y,W)]\lra
[(P,W),(X,W)]$ is injective.
\end{proposition}

\begin{proof}
Suppose given morphisms $g_0, g_1:P \longrightarrow Y$ such that
$wg_0$ and $wg_1$ are filtered homotopic. We want to prove that
$g_0$ and $g_1$ are filtered homotopic. If we define $g=g_1-g_0$,
it is enough to prove that $g$ is filtered null-homotopic. By
hypothesis, $wg$ is filtered null-homotopic, that is, there exists
a filtered homotopy $h:0\stackrel{\sim}\ra wg$. Hence the pair
$(g,h)$ defines a filtered morphism from $P$ to the path complex
$L(w)$.

We recall that   the path complex $L(f)$ of a filtered morphism
$f:B\lra A$ is defined to be the filtered complex  $L(f)= B\oplus
A[-1]$ with the differential
$$
\partial^{L(f)}=\begin{pmatrix}
\partial^B&0\\f&-\partial^A
\end{pmatrix}.
$$
 Hence, a  morphism $\alpha=\begin{pmatrix}
\alpha_0\\\alpha_1
\end{pmatrix}: X\lra L(f)$ is defined by a   morphism $\alpha_0:X\lra B$
and a filtered homotopy $\alpha_1:X\lra A[-1]$,
$\alpha_1:0\stackrel{\sim}\ra f\alpha_0$.

 On
the other hand, $\id_Y\times w:L(\id_Y)\lra L(w)$ is a filtered
quasi-isomorphism. By Proposition \ref{i-filtrado}, there exists a
lifting $(g',h'):P\lra L(\id_Y)$, up to a filtered homotopy, of
$(g,h)$.
$$
\xymatrix{ &&L(\id_Y)\ar[dd]^ {\id_Y\times w}\\ \\
P\ar@{.>}[rruu]^{(g',h')}\ar[rr]^{(g,h)}&&L(w) }
$$
That is, $(\id_Y\times w)\circ (g',h')=(g',wh')$ is filtered
homotopic to $(g,h)$. In particular, $g'$ is filtered homotopic to
$g$. In addition, $h'$ is a filtered homotopy from $0$ to $g'$,
hence $g$ is filtered null-homotopic.
\end{proof}

\begin{proposition}\label{fil2}
For each filtered complex $X$  there exists a filtered
quasi-isomorphism  $\varepsilon:P\lra X$, where $P\in\mathcal P$.
\end{proposition}

\begin{proof}
We prove, by induction,  that, for each $p$, there exists a
filtered quasi-isomorphism,  $\varepsilon_p:P_p\lra W_pX$  where
$P_p\in\mathcal P$, together with a monomorphism $j_p:P_{p-1}\lra
P_p$ such that $W_qP_{p-1}=j_p^{-1}(W_qP_p)$ for each $q$,
$W_{p}P_p=P_p$, and $\varepsilon_p$ is an extension of
$\varepsilon_{p-1}$. So the result follows from the regularity of
the filtration $W$ on $X$.

As $W$ is bounded below, we can take $P_p=0$, $\varepsilon_p=0$,
for $p\ll 0$. Assume that there is a filtered quasi-isomorphism
$\varepsilon_{p-1}: P_{p-1}\lra W_{p-1}X$, such that $ P_{p-1}\in
\mathcal P $, and $W_{p-1}P_{p-1}=P_{p-1}$. We want to extend this
model to a model of $W_pX$.

By composing $\varepsilon_{p-1}$ with the inclusion $\iota_p:
W_{p-1}X\lra W_pX$, we get a filtered morphism $\rho_{p-1}:
P_{p-1}\lra W_{p}X$. If $L\rho_{p-1}$ denotes its path complex
with the induced filtration, the complex $W_{q}L\rho_{p-1}$ is
acyclic for each $q<p$, and $W_pL\rho_{p-1}=L\rho_{p-1}$. Let
$s:G_p\lra L(\rho_{p-1})$ be a quasi-isomorphism, where $G_p$ is a
bounded below complex which is projective in each degree. Endow
$G_p$ with the filtration $W$ defined by  $W_{p-1}G_p=0$, and
$W_pG_p=G_p$. Then $s$ is a filtered quasi-isomorphism.

 Let $\xi=\pi\circ s$, and consider the following commutative diagram of filtered
complexes
$$
\begin{CD}
G_p @>\xi>> P_{p-1} @. \\
@VsVV @V\rm{id}VV @. \\
L\rho_{p-1} @>\pi>> P_{p-1} @>>> W_pX
\end{CD}
$$

Since  $0\stackrel{\sim}\ra u\id_{P_{p-1}} \xi= u\pi s$, there
exists  a filtered morphism $\nu : C\xi\lra W_pX$ which completes
the previous diagram to a commutative diagram
$$
\xymatrix{ G_p ^s\ar[dd]^s\ar[r]^\xi&
P_{p-1}\ar[dd]^{\id}\ar[r]^{j_p} &C\xi \ar[d]^\nu
\\
&&C(\pi)=Cyl(\rho_{p-1})\ar[d]^\beta\\
L\rho_{p-1}\ar[r]^{\pi}& P_{p-1}\ar[ur]^u \ar[r]^{\rho_{p-1}}&
W_pX  }
$$
where $\beta$ is a filtered homotopy equivalence (see \cite{GMa}
Lemma III.3.3). As $s$ is a filtered quasi-isomorphism, so is
$\nu$, hence we may take $P_p=C\xi$ with the induced filtration,
and $\varepsilon_p=\beta\nu$. It is clear that $W_pP_p=P_p$, and
the inclusion $j_p:P_{p-1}\lra P_{p}$ satisfies
$W_qP_{p-1}=j_p^{-1}(W_qP_p)$, for each $q$.
\end{proof}

\begin{remark}
The  equivalence of categories $\mathcal P/\sim\; \lra
\mathbf{DF}_+(\A)$ of Theorem \ref{CEfiltrados} is a well-known
result of Illusie (see \cite{I} Cor. (V.1.4.7)).
\end{remark}

\subsection{Filtered Algebras}

 In this section we review, using the formalism of
Cartan-Eilenberg categories, the homotopy theory of filtered cdg
algebras ($(R,r)$-algebras), which Halperin and Tanr\'{e} developed in
\cite{HT} by perturbation methods.

Let $R$ be a commutative ring such that $\mathbb Q\subset R$, and
$r\ge 0$ an integer. Let $\mathbf{F_rAlg}(R)$ be the category
of $(R,r)$-algebras in the sense of Halperin-Tanr\'{e} \cite{HT}.

An $(R,r)$-{\it extension\/} is a morphism of filtered cdg
algebras of the form

$$
A \longrightarrow A\widehat{\otimes} \Lambda Y \ , \qquad a
\mapsto a \otimes 1  \ ,
$$

where $\,\widehat{}\,$ means the completion of the tensor product,
$Y$ is a projective $R$-module and the morphism satisfies a
certain kind of nilpotence condition (see \emph{op.cit.}
Definition 2.2). These morphisms play the role of cofibrations in
a Quillen model category, because of the lifting property
\cite{HT}, Theorem 5.1.

An $(R,r)$-{\it quasi-isomorphism\/} is a morphism $\phi : A
\longrightarrow A'$ of $(R,r)$-algebras such that $E_{r+1}(\phi)$
is an isomorphism (here $E_{r+1}$ means the $r+1$ stage of the
associated spectral sequence). Let $\W$ denote the class of
$(R,r)$-quasi-isomorphisms.

 Given an $(R,r)$-algebra $C$, let us consider
the $(R,r)$-algebra $C\widehat{\otimes} C$. Denote by

$$
\lambda_0, \lambda_1 : C \longrightarrow C\widehat{\otimes} C \ ,
\qquad \lambda_0(c) = c\otimes 1  ,\ \lambda_1 (c) = 1 \otimes c
$$

the natural inclusions. The product

$$
\mu : C\widehat{\otimes} C \longrightarrow C
$$

defines a morphism of $(R,r)$-algebras such that $\mu \lambda_i =
\id_C , i=0,1$. A {\it cylinder object\/} for $C$ is a
factorization of $\mu$

$$
\xymatrix{ {C\widehat{\otimes} C } \ar[r]^-{\iota} &
{(C\widehat{\otimes} C)\widehat{\otimes} \Lambda X } \ar[r]^-{m} &
C }
$$

where $\iota$ is an $(R,r)$-extension and $m$ an
$(R,r)$-quasi-isomorphism. This kind of factorization exists
because of \emph{op.cit.}, Theorem 4.2. Put $\Cyl (C) =
(C\widehat{\otimes} C)\widehat{\otimes} \Lambda X$. We will also
denote by $\lambda_i$ the compositions $\iota\lambda_i$.

Let $f_0, f_1 : C \longrightarrow B$ be two morphisms of $
\mathbf{F_rAlg}(R)$. A {\it left homotopy\/} from $f_0$ to $f_1$
is a morphism $H : \Cyl (C) \longrightarrow B $ of
$(R,r)$-algebras such that $H\lambda_i = f_i \ , i=0,1$. Let $f_0
\sim f_1 $ denote the equivalence relation transitively generated
by the left homotopy.

Let us show that $\sim$ is  compatible with composition: the
implication $f_0 \sim f_1  \Longrightarrow \phi f_0 \sim \phi f_1
$ is always true for left homotopies. Let $f'_0, f'_1 : C'
\longrightarrow B'$ and $\psi : C \longrightarrow C'$ be morphisms
of $(R,r)$-algebras, such that $f'_0 \sim f'_1$. We may assume
that there is a left homotopy $H': \Cyl (C') \longrightarrow B'$
such that $H'\lambda'_i = f'_i \ , i = 0,1$. Then, consider the
following solid commutative diagram:

$$
\xymatrix{ {(C\widehat{\otimes} C)} \ar[r]^{\psi\otimes \psi}
\ar[d]^{\iota} & {(C'\widehat{\otimes} C') } \ar[r]^{\iota'} &
{\Cyl (C') }
\ar[d]^{m'}   \\
{\Cyl (C)} \ar[r]^m \ar@{-->}[rru]^h  & C \ar[r]^{\psi}  & {C'} }
$$

Because of the lifting theorem, \cite{HT}, Theorem 5.1, there
exists a morphism $h: \Cyl (C) \longrightarrow \Cyl (C')$ such
that $h\iota = \iota' (\psi\otimes\psi)$ and $m'h = \psi m$. Hence
$H= H'h$ is a left homotopy from $f'_0\psi$ to $f'_1 \psi$.
(Remark: \cite{HT} defines $\sim$ only when $A \longrightarrow C$
is an $(R,r)$-extension. In this case, left homotopy is already a
congruence by \emph{op.cit.} Proposition 6.3 and 6.5.)

Let $\Scal$ be the class of homotopy equivalences with respect to
$\sim$.

\begin{lemma}
The equivalence relation $\sim$ and $\Scal$ are compatible.
\end{lemma}

\begin{proof} Because of Example \ref{cilindro}, it is enough to show that
$m: \Cyl (C) \longrightarrow C$ is in $\Scal$. We obviously have
that $m\lambda_0 = \id_C$. Define $H: \Cyl (\Cyl (C)) = (\Cyl (C)
\widehat{\otimes} \Cyl (C)) \widehat{\otimes} \Lambda Y
\longrightarrow \Cyl (C)$ by $H(c\otimes d \otimes y) =
\lambda_0m(c) \cdot d \cdot \varepsilon (y)$, where $\varepsilon :
\Lambda Y \longrightarrow R$ is the augmentation $\varepsilon (R)
= 1, \varepsilon_{\vert Y} = 0$. Then $H$ is a homotopy from
$\lambda_0 m$ to $\id_{\Cyl (C)}$.
\end{proof}

So we have an isomorphism of categories
$\mathbf{F_rAlg}(R))/\!\sim \,\, \cong  \mathbf{F_rAlg}(R)
[\Scal^{-1}]$.

 Since, by construction, $m: \Cyl (C) \longrightarrow
C$ is in $\W$, we have $\Scal \subset \overline{\W}$. So
$(\mathbf{F_rAlg}(R), \Scal, \W )$ is a category with strong and
weak equivalences. Let $\C_{r,HT}$ be the full subcategory of $
\mathbf{F_rAlg}(R)$ of $(R,r)$-extensions of $R$.

\begin{theorem}\label{AFCE} $( \mathbf{F_rAlg}(R) , \Scal ,\W
)$ is a left Cartan-Eilenberg category and $\C_{r,HT}$ is a
subcategory of cofibrant models of $ \mathbf{F_rAlg}(R)$.
\end{theorem}

\begin{proof} Let us check the hypothesis of Theorem \ref{F}: if
$M$ is an $(R,r)$-extension of $R$ and $E \longrightarrow B$ a
quasi-isomorphism of $ \mathbf{F_rAlg}(R) $, the induced map
$[M,E] \longrightarrow [M,B]$ is a bijection by \cite{HT},
Application 7.7. So $(R,r)$-extensions are cofibrant objects. The
existence of enough cofibrant objects is guaranteed by Theorem
4.2. op.cit.
\end{proof}

 Halperin-Tanr\'{e} also define a notion of minimal
$(R,r)$-algebras (\emph{op.cit.}, Definition 8.3) and when $R =
\mk$ is a field of characteristic zero, they prove the existence
of minimal models for $(\mk,r)$-algebras $B$ such that
$E_{r+1}(B)$ is concentrated in non-negative degrees and
$H^0(E_r(B)) = \mk$. Let $ \mathbf{F_rAlg}(\mk)_{0}$ denote the
full subcategory of $ \mathbf{F_rAlg}(\mk) $ with  objects the
$(\mk,r)$-algebras $B$ such that $E_{r+1}(B)$ is concentrated in
non-negative degrees and $H^0(E_r(B)) = \mk$,  and let $\M_{r,HT}$
denote  the full subcategory of minimal $(\mk,r)$-algebras. We can
sum up their results in the following theorem.

\begin{theorem}\label{AFS} Let  $\mk$ be a field of characteristic zero.
$( \mathbf{F_rAlg}(\mk)_{0} , \Scal, \W)$ is a left Sullivan
category and $\M_{r,HT}$ is the subcategory of minimal objects of
$\mathbf{F_rAlg}(\mk)_{0}$.
\end{theorem}

\section{Cartan-Eilenberg categories defined by a cotriple}

In section 3 we have proved, under suitable hypotheses, that some
subcategories of the functor category  $\Cat(\C,\D)$ are
Cartan-Eilenberg categories, and as a consequence we shaw that the
derived functor of an additive functor $K$ is a cofibrant model of
$K$. In this section we prove that the whole category
$\Cat(\C,\D)$ is a Cartan-Eilenberg category if
 $\C$ has a cotriple  and $\D$ is a category of chain complexes.
The cofibrant model of a functor $K$ with respect to this
structure is the non-additive derived functor of $K$ as introduced
by Barr-Beck (\cite{BB1}).

\subsection{Categories of chain complexes and cotriples}\label{aciclicosBarr}

Let $\A$ be an additive category and denote by $\mathbf C_{\ge
0}(\A)$ the category of non-negative chain complexes of $\A$. In
this section  we will consider   as strong equivalences in
$\mathbf C_{\ge 0}(\A)$  classes of summable morphisms as
introduced in the following definition.

\begin{definition}\label{summable} Let $\A$ be an additive category.
A class $\mathcal S$  of morphisms of $\mathbf C_{\ge 0}(\A)$ is
called a class of \emph{summable} morphisms if it satisfies the
following  properties:
\begin{itemize}
\item[(i)] $\mathcal S$ is saturated.
\item[(ii)] The homotopy equivalences are in $\mathcal S$.
\item[(iii)] Let $f:C_{\ast \ast }\lra D_{\ast \ast}$ be a morphism of  first
quadrant double  complexes. If $f_n:C_{\ast n }\lra D_{\ast n}$ is
in $\mathcal S$ for all $n\ge 0$, then $\Tot f: \Tot C_{\ast \ast
}\lra \Tot D_{\ast \ast}$ is in $\mathcal S$.
\end{itemize}
\end{definition}

\begin{nada}
For example, the class of homotopy equivalences, which will be
denoted by $\mathcal S_h$, is a class of summable morphisms. Also,
if $\A$ is an abelian category,   the class of quasi-isomorphisms
is a class of summable morphisms (cf. \cite{B}, Chap. 5).
\end{nada}

\begin{nada} Let $\A$ be an additive category,  and let
$$\mathbf G=(G:\A\ra \A,\varepsilon:G\Rightarrow \id_\A,
\delta:G\Rightarrow G^2)$$ be a cotriple  on  $\A$.

We recall that the cotriple  $\mathbf G$ is called \emph{additive
} if the functor $G$ is additive, in such case, it induces an
additive cotriple on $\mathbf C_{\ge 0}(\A)$ which we also denote
by $\mathbf G$.

Let  $\mathcal S$ be a class of summable morphisms of $\mathbf
C_{\ge 0}(\A)$, and $\mathbf G$ an additive cotriple on $\A$. We
say that $\mathbf G$ and $\mathcal S$   are \emph{compatible}  if
the extension of $G$  to the category of complexes  $G: \mathbf
C_{\ge 0}(\A)\lra \mathbf C_{\ge 0}(\A)$ satisfies $G(\mathcal
S)\subset {\mathcal S}$. In this case, taking $\W=G^{-1}({\mathcal
S})$, $(\mathbf C_{\ge 0}(\A),\mathcal S,\W)$ is a category with
strong and weak equivalences.

For example, the class of homotopy equivalences  $\mathcal S_h$ in
$\mathbf C_{\ge 0}(\A)$ is compatible with any additive cotriple
$\mathbf G$ on $\A$, thus, taking  $\W_h= G^{-1}(\mathcal S_h)$,
$(\mathbf C_{\ge 0}(\A),\mathcal S_h,\W_h)$ is a category with
strong and weak equivalences.

\end{nada}

\begin{nada}\label{aciclicdos}
Let  $\mathbf G = (G,\varepsilon, \delta )$ be an additive
cotriple defined on the category $\A$, and by extension on
$\mathbf C_{\ge 0}(\A)$.

The simplicial standard construction associated to the cotriple
$\mathbf G$ on $\mathbf C_{\ge0}(\A)$ defines, for each object $K$
in $ \mathbf C_{\ge 0}(\A)$,
 an augmented  simplicial object  $\varepsilon:B_\bullet (K)\lra K$  in $ \mathbf C_{\ge 0}(\A)$ such that
$B_n(K) = G^{n+1}(K)$, (\cite{G}, App.,  see also \cite{ML}).
Hence, there is a naturally defined double complex $B_{*}(K)$
associated to $B_\bullet (K)$, with total complex $B(K)= \Tot
B_*(K)$. This construction defines a functor
$$
B:  \mathbf C_{\ge 0}(\A)\lra \mathbf C_{\ge 0}(\A),
$$
with a natural transformation $\varepsilon : B\Rightarrow 1$.

\end{nada}

\begin{theorem}\label{modelsaciclics1}
Let   $\A$ be an additive category, $\mathbf G$ an additive
cotriple
 on $\A$, and  $\mathcal S$ a class of summable
morphisms in  $\mathbf C_{\ge 0}(\A)$  compatible with $\mathbf
G$.
 Then, with the previous notation,
\begin{enumerate}[(1)]
\item $(B,\varepsilon )$ is a left resolvent functor for
$ (\mathbf C_{\ge 0}(\A), \mathcal{S}, \W)$;

\item $( \mathbf C_{\ge 0}(\A), \mathcal{S}, \W)$ is a left
Cartan-Eilenberg category; and

\item an object $K$ of  $\mathbf C_{\ge 0}(\A)$
is cofibrant if and only if $\varepsilon_K:B(K)\lra K$ is in
$\mathcal S$.

\end{enumerate}

\end{theorem}
\begin{proof}
Let us verify the hypotheses of Theorem \ref{MA}. Firstly, if
$s\in \mathcal S$, then $G(s)\in {\mathcal S}$ by hypothesis, and
it follows, inductively, that $G^i(s)\in {\mathcal S}$ for any
$i\geq 0$. By  Definition \ref{summable}(iii), we deduce that
$B(s)=\Tot B_{*}(s)\in {\mathcal S}$. Therefore $B(\mathcal
S)\subset \mathcal S$.

Next, let $K$ be a chain complex of $\A$. For any $i>0$, the
augmented simplicial objects $\varepsilon_{G^iK}: B_\bullet
G^iK\lra G^iK$ and $G^i(\varepsilon_K):G^iB_\bullet K\lra G^iK$
have a contraction induced by the morphism $\delta:G\lra G^2$,
hence $\varepsilon_{G^iK},G^i(\varepsilon_K)\in \mathcal S$, by
\ref{summable}(ii). Therefore, by \ref{summable}(iii) and the
additivity of $G$,  the morphisms $$\varepsilon_{G^iK}: B G^iK\lra
G^iK\quad \text{and }\quad G^i(\varepsilon_K):G^i{\rm{Tot}}B_*
K\cong {\rm{Tot}}G^i B_* K\lra G^iK$$  are in $\mathcal S$, for
each $i>0$. Applying again \ref{summable}(iii)  we obtain that
$B(\varepsilon_K)$ and $\varepsilon_{BK}$ are in $\mathcal S$.
Therefore, $(B,\varepsilon)$ is a left resolvent functor for $
(\mathbf C_{\ge 0}(\A), \mathcal{S}, B^{-1}(\mathcal S))$, by
Theorem \ref{MA}.

Finally, let us check that $\W=B^{-1}({\mathcal S})$, that is
$G^{-1}(\mathcal S)=B^{-1}(\mathcal S)$. Indeed, let $w:K\lra L$
be a morphism of functors. If $w\in \W=G^{-1}({\mathcal S})$, we
have $G^i(w)\in {\mathcal S}$ for each $i>0$, therefore, applying
once again \ref{summable}(iii), we obtain $B(w)\in {\mathcal S}$.
Conversely, if $B(w)\in{\mathcal S}$, since $BGw=GBw$, we have
$BGw\in {\mathcal S}$, and from the commutativity of the diagram
$$\xymatrix{BGK\ar[d]_{\varepsilon_{GK}}\ar[r]^{BGw}&BGL
\ar[d]^{\varepsilon_{GL}}\\
GK\ar[r]^{Gw}&GL}$$ it follows that $Gw\in {\mathcal S}$, because
$\varepsilon_{GK},\varepsilon_{GL}\in \mathcal S$, by
\ref{summable} (ii), and  $\mathcal S$ is saturated, by
\ref{summable} (i). Hence $w\in \W$.
\end{proof}

In order to recognise cofibrant objects in $(\mathbf C_{\ge
0}(\A),\mathcal S,\W)$
 the following
criterion   will be useful.

\begin{proposition}\label{critcof}
Let  $\A$ be an additive category, $\mathbf G$ an additive
cotriple on $\A$, and $\mathcal S$ a class of summable morphisms
in $\mathbf C_{\ge 0}(\A)$ compatible with  $\mathbf G$. Then,

\begin{enumerate}\item for each object $K$ of $\mathbf C_{\ge
0}(\A)$, $GK$ is cofibrant;

\item if $K$ is an object of $\mathbf C_{\ge 0}(\A)$  such that
$K_{n}$ is cofibrant for each $n\ge 0$, then  $ K$ is cofibrant
(in \cite{B}  one such complex is called
$\varepsilon$-presentable); and

\item if $K$ is an object of $\mathbf C_{\ge 0}(\A)$ such that
$\varepsilon_{K_n}: G(K_n)\lra K_n$ has a section, that is to say,
there are morphisms $\theta_n:K_n\lra G(K_n)$ such that
$\varepsilon_{K_n}\theta_n =\id_{K_n}$, for  $n\ge 0$, then, $K$
is cofibrant (in \cite{BB1}  one such complex is called
$G$-representable).

\end{enumerate}
\end{proposition}

\begin{proof}
(1) The augmented simplicial complex $\varepsilon_{GK}:B_\bullet
GK\lra GK$ is contractible, because the morphism $\delta_K:GK\lra
G^2K$ induces a contraction. Hence, by \ref{summable}(ii),
$\varepsilon_{GK}\in \mathcal S$, so $GK$ is cofibrant, by Theorem
\ref{modelsaciclics1}(3).

(2) Suppose  $K_n$  cofibrant, for each $n\ge 0$. Then
$\varepsilon_{K_n}:B(K_{n})\lra K_{n}\in \mathcal S$, by Theorem
\ref{modelsaciclics1}(3). Therefore $\varepsilon_K:BK\lra K\in
\mathcal S$, by \ref{summable}(iii), hence $K$ is cofibrant, by
Theorem \ref{modelsaciclics1}(3) again.

(3) Each $G(K_{n})$ is cofibrant, by (1), and $K_n$ is a retract
of $G(K_n)$,  then, by Proposition \ref{retracto}, $K_n$ is
cofibrant. Hence, by (2), $K$ is cofibrant.
\end{proof}

\subsection{Functor categories and cotriples}

\begin{nada}\label{ejemploaciclicas}
Given a category $\X$ and an additive category $\A$,  the functor
category $\Cat(\X, \A)$ is also additive, so we can consider
classes of summable morphisms in
$$
\mathbf C_{\ge 0}\Cat(\X, \A)\cong\Cat(\X, \mathbf C_{\ge 0}(\A)
).
$$
Besides the class of natural homotopy equivalences $\mathcal S_h$,
 we will consider point-wise defined classes. Take $\Sigma$ a class of summable morphisms in
$\mathbf C_{\ge 0}(\A)$ and define a class of morphisms $\mathcal
S_\Sigma$ of $\Cat(\X, \mathbf C_{\ge 0}(\A) )$ by
$$
\mathcal S_\Sigma =\{\ f\ ;\ \ f(X)\in\Sigma,\ \ \forall X\in
Ob\X\ \}.
$$
Then $\mathcal S_\Sigma $ is a class of summable morphisms. We
shall say that $\mathcal S_\Sigma $ is the class of summable
morphisms in $\Cat(\X, \mathbf C_{\ge 0}(\A) )$ defined
\emph{point-wise} from $\Sigma$.

For example, if $\Sigma$ is the class of homotopy equivalences in
$\mathbf C_{\ge 0}(\A)$, we say that $\mathcal S_\Sigma $ is the
class of \emph{point-wise homotopy equivalences}, and we denote it
by $\mathcal S_{ph}$. Observe that in contrast to the case of
natural homotopy equivalences $\mathcal S_h$ in $\Cat(\X, \mathbf
C_{\ge 0}(\A) )$, the point-wise homotopy equivalences have
homotopy inverses over each object $X$ of $\X$, but these homotopy
inverses are not required to be natural.  So, generally speaking,
the inclusion $\mathcal S_{h}\subset \mathcal S_{ph}$ is strict.

\end{nada}\begin{nada} If $\mathbf{G}$ is a cotriple in $\X$,
it naturally defines an additive cotriple on the functor category
$\Cat(\X, \mathbf C_{\ge 0}(\A) )$  by sending $K$ to $K\circ G$,
with the evident extensions of the transformations $\varepsilon ,
\delta$. We also denote this cotriple by $\mathbf{G}$.

If $\mathcal S$ is   a class of point-wise defined morphisms , then
$\mathcal S$ is compatible with each cotriple $G$ on $\Cat(\X,
\mathbf C_{\ge 0}(\A) )$ induced by a cotriple on $\X$.

\end{nada}

From Theorem \ref{modelsaciclics1}, taking $\W$ as above, we
obtain the following result.

\begin{theorem}\label{MACatAd} Let $\X$ be a category and $\A$ an additive
category.  Let  $\mathbf G$ be  an additive  cotriple on
$\Cat(\X,\A)$, and  $\mathcal S$  a class  of summable morphisms
in $\Cat(\X, \mathbf  C_{\ge 0}(\A))$ compatible with $\mathbf G$.
Then,
\begin{enumerate}

\item $(B,\varepsilon)$ is a left resolvent functor for
$(\Cat(\X, \mathbf  C_{\ge 0}(\A)),\mathcal S,\W)$;

\item $(\Cat(\X, \mathbf  C_{\ge 0}(\A)),\mathcal S,\W)$ is a left
Cartan-Eilenberg category; and

\item an object $K$ of $\Cat(\X, \mathbf  C_{\ge 0}(\A))$ is
cofibrant if and only if $\varepsilon_K:BK\lra K$ is in $\mathcal
S$.
\end{enumerate}

\end{theorem}

\begin{nada} Let $\X$ be a category with arbitrary sums.
We recall that, associated to each set $\M$  of objets of $\X$
(called ``models"), there is defined a cotriple $\mathbf G$ on
$\X$ (see for example \cite{BB2},  (10.1)). The functor $G$ is
given by the formula
$$
G(X)= \bigsqcup_{ (M,\ f)\in \X/X} M_f ,
$$
where, if $f:M\lra X$ is an object of $\X$ over $X$,  $M_f$ denotes
a copy of $M$ indexed by $f$. Denote by $\langle f\rangle :M\lra
G(X)$ the canonical inclusion  into the sum corresponding to the
summand $M_f$. If $a:X\lra Y$ is a morphism, $G(a):G(X)\lra G(Y)$
is defined in such a way that $G(a)\circ \langle f\rangle =\langle
af\rangle $, for each $f:M\lra X$. The counit
$\varepsilon:G\Rightarrow 1$ is defined by $\varepsilon_{X}\circ
\langle f\rangle =f$, and comultiplication $\delta:G\Rightarrow
G^2$, is defined by $\delta_{X}\circ \langle f\rangle =\langle
\langle f\rangle \rangle $.

\end{nada}
\begin{nada}\label{cotriplefunctores}
In the same way, for a general category $\X$ with a set $\M$ of
objects, if the additive category $\A$ has arbitrary sums, there
is a variant of the model-induced cotriple given as follows. The
cotriple  $\mathbf G$ in $\Cat(\X,\A)$ is defined by
$$
(GK)(X)=\bigoplus_{(M,\ f)\in\X/X}K(M_f),\;
$$
with counit $\varepsilon:G\Rightarrow 1$ defined by
$\varepsilon_{K,X} \circ\langle f\rangle =K(f)$, and
comultiplication $\delta:G\Rightarrow G^2$, defined by
$\delta_{K,X}\circ \langle f\rangle =\langle \langle f\rangle
\rangle $. This cotriple  is additive.

\end{nada}

\begin{example}\label{notationG}
In the  original formulation of the  Beck homology (see
\cite{BB2}), one considers

\begin{itemize}
\item[(a)] a cotriple $\mathbf G$ defined on the category $\X$,
\item[(b)] an abelian category $\A$,
\item[(c)] a class of acyclic morphisms in $\Cat(\X,\mathbf C_+(\A))$, and
\item[(d)] a functor $F:\X\lra \A$.
\end{itemize}

Then, the homology of $X$ with coefficients in $F$ is defined as
$H_*(X,F)_\mathbf{G}=H_*((BF)(X))$, that is, the homology of the
cofibrant model of $F$.

In this case, the category $\A$ is abelian and the cotriple on
$\Cat(\X,\A)$  is induced by a cotriple on $\X$.
\end{example}

\begin{examples}\label{ejemplotop}

(1) Barr-Beck proved that the singular homology functors with
integer coefficients $H_* = \left\{H_n \right\}_{n= 0,1,\dots}$
are the derived functors of the $0$-th singular homology functor
$H_0$. In this first example we give a version of this result at
the chain level: we prove that the functor of singular chains
$S_*$ is a cofibrant model for the functor $H_0$ in the category
of chain complex valued functors on topological spaces with a
convenient Cartan-Eilenberg structure.

Let $\X=\Top$ be the category of topological spaces and consider
the  cotriple $\mathbf G$ on $\Top$ defined by the set
$\{\Delta^n; n\in \mathbb N\}$,
$$
G(X)=\bigsqcup_{ (\Delta^n,\ \sigma)\in \mathbf{Top}/X}
\Delta_\sigma^n.
$$
We consider  on the category of abelian groups valued functors
$\Cat(\Top,\mathbb Z\mathbf{-mod})$ the cotriple induced by
$\mathbf{G}$.

Take $\mathcal S_h$ the class of natural homotopy equivalences in
$\Cat(\Top,\mathbf C_{\ge 0}(\mathbb Z))$ and
$\W_{h}=G^{-1}(\mathcal S_h)$. From Theorem \ref{modelsaciclics1}
we obtain that  $(\Cat(\Top, \mathbf C_{\ge 0}(\mathbb Z)),
\mathcal S_h, \W_h)$ is a left Cartan-Eilenberg category.

Let $S_\ast:\Top\lra \mathbf C_{\ge 0}(\mathbb Z)$  be the functor
of singular chains with integer coefficients, and $\tau:
S_\ast\lra H_0$ the natural augmentation.

Let us see  that $S_*$ is cofibrant.  Let $\theta_n:S_n\lra
S_n\circ G$ be the natural transformation which, for each
topological space $X$, sends a singular simplex
$\sigma:\Delta^n\lra X$ to $\theta_n(\sigma) =
\langle\sigma\rangle$. It is clear that $\varepsilon_{S_n}
\theta_n =\id_{S_n}$, so $S_\ast$ is cofibrant, by Proposition
\ref{critcof}(iii).

On the other hand, the morphism $\tau: S_\ast\lra H_0(-,\mathbb
Z)$ is in $\W_h$. In fact, for each $n\ge 0$, take a homotopy
inverse of $\tau_{\Delta^n}$, $\lambda_n: H_0(\Delta^n,\mathbb
Z)\lra S_*(\Delta^n)$. Then, for each topological space $X$,
$$\lambda_X:=\bigoplus_{(\Delta^n,\ \sigma)\in \mathbf{Top}/X}
(\lambda_{n},\sigma):H_0(GX,\mathbb Z)\lra S_*(GX)$$ defines a
natural morphism $\lambda:  H_0(-,\mathbb Z)\lra S_*$ which is a
homotopy inverse of $\tau$.

Hence, $S_\ast$ is a cofibrant model for $H_0(-,\mathbb Z)$ in
$(\Cat(\Top, \mathbf C_{\ge 0}(\mathbb Z)), \mathcal S_h, \W_h)$.

Notice that, if $\mathcal S$ denotes the homotopy equivalences and
$\W$ the weak homotopy equivalences in $\Top$, then
$(\Top,\mathcal S,\W)$ is a left Cartan-Eilenberg category. If we
consider in $\mathbf C_{\ge 0}(\mathbb Z)$ the class $\E$ of the
quasi-isomorphisms, the category of functors $\Cat_{\mathcal
S,\E}(\Top,\mathbf C_{\ge 0}(\mathbb Z))$ (see \ref{notacion3.3}
for the notation) has a structure of left Cartan-Eilenberg
category for which  the functor $H_0$ is a cofibrant object and
$S_*\lra H_0$ is not a weak equivalence.

{(2)} The next example is a variation for differentiable manifolds
of the previous one.

Let $\X=\Diff$ be the category of differentiable manifolds with
corners. Consider the additive cotriple $\mathbf G^\infty$ defined
on $\Cat(\Diff, {\mathbb Z}\mathbf{-mod})$ by the set $\{\Delta^n;
n\in \mathbb N\}$,
$$
G^\infty(K)(X)=\bigoplus_{(\Delta^n,\ \sigma)\in \mathbf{Diff}/X}
K(\Delta^n,\sigma).
$$
By Theorem \ref{modelsaciclics1}, $\Cat(\Diff, \mathbf C_{\ge
0}(\mathbb Z)), \mathcal S_h, \W_{h})$ is a left Cartan-Eilenberg
category.

Denote by $S_\ast^\infty :\Diff\lra \mathbf C_{\ge 0}(\mathbb Z)$
the functor of \emph{differentiable} singular chains. Reasoning as
in the topological case, it follows that $S^\infty_\ast$ is a
cofibrant model of $H_0(-,\mathbb Z)$ in the left Cartan-Eilenberg
category $(\Cat(\Diff, \mathbf C_{\ge 0}(\mathbb Z)), \mathcal
S_h, \W_h)$.

These two previous examples permit us to give an interpretation of
a well-known theorem of Eilenberg for the singular complex of a
differentiable manifold (see \cite{E} and its extension to
differentiable manifolds with corners in \cite{Hu}). By
Eilenberg's theorem the natural transformation $S^\infty_\ast\lra
S_\ast$ is a point-wise homotopy equivalence in $\Cat(\Diff,
\pChains{\mathbb Z})$, hence $S_\ast$ is  a cofibrant model of
$H_0(-,\mathbb Z)$ in
 $(\Cat(\Diff, \mathbf C_{\ge 0}(\mathbb Z)), \mathcal S_{ph},
\W_{ph})$. However, $S^\infty_\ast$ and $S_\ast$ are not naturally
homotopy equivalent functors  in $\Cat(\Diff,\mathbf C_{\ge
0}(\mathbb Z))$ (see \cite{GNPR3}), so $S_\ast$ is not a cofibrant
model of $H_0(-,\mathbb Z)$ in $(\Cat(\Diff, \mathbf C_{\ge
0}(\mathbb Z)), \mathcal S_h, \W_h)$.

Observe that the
 Cartan-Eilenberg category $(\Cat(\Diff,
\mathbf C_{\ge 0}(\mathbb Z)), \mathcal S_{ph}, \W_{ph})$  does not
come from a Quillen model category, since the morphisms in the
class $\mathcal S_{ph}$ do not have, in general, a homotopic
inverse.

\end{examples}

\begin{nada} If, in  Theorem \ref{MACatAd},
 the cotriple $\mathbf G$ is induced by a cotriple on $\X$, we
can prove that the natural transformations from a cofibrant
functor $K$ to any other functor $L$  are determined by its
restriction to the ``models", as stated in the following theorem.

 Let $\X$ be a category with a
cotriple $\mathbf G$,  let $\A$ be an additive category, and
$\mathcal S$  a class of summable morphisms in $\Cat(\X, \mathbf
C_{\ge 0}(\A))$ compatible with the additive cotriple induced by
$\mathbf G$.

 Denote by  $\M$ the full subcategory of
$\X$ with objects $GX$, for $X\in Ob\X$ and by
$$\rho:\Cat(\X,\mathbf C_{\ge 0}(\A))\lra \Cat(\M,\mathbf C_{\ge
0}(\A))$$ the restriction functor, $\rho(K)=K_{|\M}$.

Since $G$ sends objects in $\M$ to $\M$, $\mathbf G$ induces a
cotriple on $\M$, and a functor $B_\M$ such that $\rho \circ
B=B_\M\circ \rho $.

Since $B_\bullet:\X\lra \Delta^{op}\X$ factors through the
inclusion $\Delta^{op}\M\lra  \Delta^{op}\X, $ the functor
$$B:\Cat(\X,\mathbf C_{\ge 0}(\A))\lra \Cat(\X,\mathbf C_{\ge
0}(\A)),\quad BK=Tot\circ \Delta^{op}K\circ B_\bullet$$  factors
through $\rho$, that is, if $$B':\Cat(\M,\mathbf C_{\ge
0}(\A))\lra \Cat(\X,\mathbf C_{\ge 0}(\A))$$ is defined by  $$ B'
K=Tot \circ\Delta^{op}K\circ  B_\bullet, $$ then $B=B'\circ \rho$.
In addition, $\rho \circ B'=B_\M$.

We say that a class  $\mathcal S_{\M}$ of morphisms in
$\Cat(\M,\mathbf C_{\ge 0}(\A))$ is \emph{adapted} to
$(\mathbf{G},\mathcal{S})$ if $\rho(\mathcal S)\subset \mathcal
S_\M$ and $B'(\mathcal S_\M)\subset \mathcal S$.  In that case,
the restriction $\rho$ induces a functor
$$
\overline{\rho}:\Cat(\X,\mathbf C_{\ge 0}(\A))[\mathcal
S^{-1}]\lra\Cat(\M,\mathbf C_{\ge 0}(\A))[\mathcal S_\M^{-1}],
$$
 and the functor $B'$
induces a functor
$$\beta':\Cat(\M,\mathbf C_{\ge 0}(\A))[\mathcal S_\M^{-1}]\lra
\Cat(\X,\mathbf C_{\ge 0}(\A))[\mathcal S^{-1}],$$such that
$\overline{\rho}\circ \beta'=\beta_\M,$ and $ \beta'\circ
\overline{\rho}=\beta$, where $\beta$ and $\beta_\M$ denote the
functors induced by $B$ and $B_\M$, respectively.

If there exists a class $\mathcal S_\M$ adapted to
$(\mathbf{G},\mathcal{S})$ we say that $\mathcal S$ is
\emph{adaptable} to $\mathbf G$.

For example, if $\mathcal S$ is the class of homotopy
equivalences, then $\mathcal S$ is adaptable to any cotriple
$\mathbf G$ on $\X$, since it is enough to take $\mathcal S_\M$ as
the class of homotopy equivalences. On the other hand, if
$\mathcal S$ is defined point-wise by a class $\Sigma$, then
$\mathcal S$ is also adaptable, taking the class $\mathcal S_\M$
 point-wise defined by  $\Sigma$.

If $K,L$ are objects in $\Cat(\X,\mathbf C_{\ge 0}(\A))$ (resp.
$\Cat(\M,\mathbf C_{\ge 0}(\A))$) we denote by $[K,L]$  the
morphisms from $K$ to $L$ in the category  $\Cat(\X,\mathbf C_{\ge
0}(\A))[\mathcal S^{-1}]$ (resp. $\Cat(\M,\mathbf C_{\ge
0}(\A))[\mathcal S_\M^{-1}]$).

\end{nada}

\begin{theorem}\label{Teorema aditivo modelos}
Let $\X$ be a category with a cotriple $\mathbf G$,  let $\A$ be
an additive category, and  $\mathcal S$  a class of summable
morphisms in $\Cat(\X, \mathbf C_{\ge 0}(\A))$ compatible with the
cotriple induced by $\mathbf G$ and $S_\M$ a class of morphisms in
$\Cat(\M, \mathbf C_{\ge 0}(\A))$ adapted to
$(\mathbf{G},\mathcal{S})$. If $K$ is a cofibrant object
 of $\Cat(\X, \mathbf C_{\ge 0}(\A))$,   the restriction map
$$ \overline{\rho}_{_{KL}}:[K,L]\lra [K_{|\M},L_{|\M}]
$$
is bijective, for each $L$.
\end{theorem}

\begin{proof}
The diagram
$$
\xymatrix{
[K,L]\ar[ddrr]^{\beta_{KL}}\ar[rr]^{\overline{\rho}_{KL}}
&&[K_{|\M},L_{|\M}]\ar[rrdd]^{\beta_{\M,K_{|\M},L_{|\M}}}\ar[dd]^{\beta'_{KL}}&&\\ \\
&&[BK,BL]\ar[rr]^{\overline{\rho}_{BK,BL}} &&[BK_{|\M},BL_{|\M}]}
$$
is commutative, since  $\beta'\circ \overline{\rho}=\beta$ and
$\overline{\rho}\circ \beta'=\beta_M$.

 By the
naturality of $\varepsilon:B\Rightarrow \id$, the following
diagram
$$
\xymatrix{[K,L]\ar[ddrr]^{\varepsilon^*_K}
\ar[rr]^{\beta_{KL}}&&[BK,BL]\ar[dd]^{\varepsilon_{L*}}
\\
\\&&[BK,L]}
$$
is commutative. Since $BK$ is cofibrant and $\varepsilon_{L}$ is a
weak equivalence, the map $ \varepsilon_{L*}$ is bijective. Since
$K$ is cofibrant, $\varepsilon_{K}$ is a strong equivalence, so
$\varepsilon_{K}^{*}$ is also bijective,
 hence
$\beta_{KL}$  is bijective. In particular, $\beta'_{KL}$ is
surjective.

On the other hand,
$\varepsilon_{K_{|\M}}={\rho}(\varepsilon_K):BK_{|\M}\lra K_{|\M}$
is in $\mathcal S_{\M}$, since $\rho(\mathcal S)\subset \mathcal
S_\M$, so $(\varepsilon_{K_{|\M}})^{*}$ is bijective. From
$(\varepsilon_{L_{|\M}})_*\circ
\beta_{\M,K_{|\M},L_{|\M}}=(\varepsilon_{K_{|\M}})^*$, we obtain
that $\beta_{\M,K_{|\M},L_{|\M}}$ is injective, so too is
$\beta'_{KL}$.

Since $\beta'_{KL}$ and $\beta_{KL} $ are  bijective maps, so too
is $\overline{\rho}_{{KL}}:[K,L]\lra [K_{|\M},L_{|\M}]$.
\end{proof}

\begin{corollary} Under the hypothesis of the previous theorem,
let $K,L$ be cofibrant   objects of $\Cat(\X, \mathbf C_{\ge
0}(\A))$. If $K_{|\M}$ and $L_{|\M}$ are isomorphic in $\Cat(\M,
\mathbf{C}_{\ge 0}(\A))[\mathcal{S}_\M^{-1}]$, then $K$ and $L$
are isomorphic in  $\Cat(\X, \mathbf{C}_{\ge
0}(\A))[\mathcal{S}^{-1}]$.
\end{corollary}

\begin{nada}
Denote the full subcategory of $\X$ defined by $\M$ by the
same letter $\M$. If $\Sigma$ is a class of morphisms of morphisms
of $ \mathbf C_{\ge 0}(\A)$ stable by arbitrary sums, the class
$\mathcal S_\Sigma$ of morphisms of $\Cat(\X, \mathbf C_{\ge
0}(\A))$ defined point-wise by $\Sigma$ is compatible with the
cotriple $\mathbf G$ associated to $\M$ by
\ref{cotriplefunctores}, and it is also adaptable to it. The class
$\mathcal S_h$ of homotopy equivalences is also adaptable to
$\mathbf G$.
\end {nada}

The  proof of Theorem \ref{Teorema aditivo modelos} works for the
following variant.

\begin{theorem} Let $\X$ be a category and $\M$ a set of objets of
$\X$. Let $\A$ be an additive category, and  $\mathcal S$  a class
of summable morphisms in $\Cat(\X, \mathbf C_{\ge 0}(\A))$
compatible with the cotriple on $\Cat(\X, \mathbf C_{\ge 0}(\A))$
associated to $\M$ and $\mathcal S_\M$ a class of morphisms in
$\Cat(\M, \mathbf C_{\ge 0}(\A))$ adapted to
$(\mathbf{G},\mathcal{S})$. If $K$ is  a cofibrant object
 of $\Cat(\X, \mathbf C_{\ge 0}(\A))$,   the map
$$ \rho_{_{KL}}:[K,L]\lra [K_{|\M},L_{|\M}]
$$
is bijective, for each $L$.
\end{theorem}

\begin{nada} The Barr-Beck's acyclic models theorem is stated in terms of acyclic
functors. We introduce the corresponding notions in our setting.
\end{nada}

\begin{definition} Let $\A$ be an abelian category, and
 $\mathcal S$  a class of summable
morphisms in  $ \mathbf C_{\ge 0}(\A)$.  If the morphisms in
$\mathcal S$ are quasi-isomorphisms, $\mathcal S$ is called a
class of \emph{acyclic} morphisms (see \cite{B}, Chap. 5, (1.1)
AC-4).

Let $\mathbf G$ be a cotriple on  $ \mathbf C_{\ge 0}(\A)$
compatible  with $\mathcal S$.
 An object $K$ of $
\mathbf  C_{\ge 0}(\A)$ is called $G$-\emph{acyclic} if the
augmentation $\tau_K: K\lra H_0K$ is a weak equivalence, that is,
$G(\tau_K)\in \mathcal S$.

\end{definition}

\begin{nada}
If $\mathcal S$ is a class of acyclic morphisms in $ \mathbf
C_{\ge 0}(\A)$, and $\phi:K\lra L$ is  morphism in $ \mathbf
C_{\ge 0}(\A)[\mathcal S^{-1}]$, then $\phi $ defines a morphism
$H_*\phi:H_*K\lra H_*L$. In particular, $H_0$ defines a functor
$$\mathbf C_{\ge 0}(\A)[\mathcal S^{-1}]\lra \A,$$
and a map $ H_0:[K,L]\lra [H_0K,H_0L]$, where $[H_0K,H_0L]$ is
simply the class of morphisms $H_0K\lra H_0L$ in $\A$.

If $K$ is a chain complex in non-negative degrees, and $L$ is a
complex concentrated in degree $0$, then the map
$$
H_0:[K,L]\lra [H_0K,H_0L]
$$
is bijective, with inverse map $ \tau_{K*}$.

\end{nada}

Now, we derive a variation of Barr-Beck's acyclic models theorem
(\cite{B}, Chap. 5, (3.1)) as a consequence of the
Cartan-Eilenberg structure of $\Cat(\X, \mathbf C_{\ge 0}(\A))$.

\begin{theorem}\label{modelosaciclicosSI}{\rm (Acyclic models theorem.)}
Let  $\X$ be a category with a cotriple $\mathbf G$, let $\A$ be
an abelian category, and  $\mathcal S$ a class of acyclic
morphisms in $\Cat(\X, \mathbf C_{\ge 0}(\A))$ compatible  with,
and adaptable to, the cotriple induced by $\mathbf G$. If $K,L$
are objects of $\Cat(\X, \mathbf C_{\ge 0}(\A))$ such that $K$ is
cofibrant and $L$ is $G$-acyclic, then the map
$$
H_0\rho_{_{KL}}:[K,L]\lra [H_0K_{|\M},H_0L_{|\M}]
$$
is bijective.
\end{theorem}

\begin{proof}
The map
$$
H_0\rho_{_{KL}}:[K,L]\lra [H_0K_{|\M},H_0L_{|\M}]
$$
factors as
$$\xymatrix{
[K,L]\ar[r]^{\tau_{L*}} &[K,H_0L]\ar[r]^{\rho}
&[K_{|\M},H_0L_{|\M}]\ar[r]^{H_0}&[H_0K_{|\M},H_0L_{|\M}] }$$

The map $\tau_{L*}$ is bijective because $K$ is cofibrant and $L$
is $G$-acyclic.  By Theorem \ref{Teorema aditivo modelos}, $\rho$
is also bijective. Finally, the map
$$H_0:[K_{|\M},H_0L_{|\M}]\lra [H_0K_{|\M},H_0L_{|\M}]
$$
is  bijective  because $K_{|\M}$ is concentrated in non-negative
degrees and $H_0L_{|\M}$ is concentrated in degree $0$.
\end{proof}

\begin{corollary} Under the hypothesis of the previous theorem,
let $K,L$ be cofibrant $G$-acyclic    objects of $\Cat(\X, \mathbf
C_{\ge 0}(\A))$. If $H_0K_{|\M}$ and $H_0L_{|\M}$ are isomorphic,
then $K$ and $L$ are isomorphic in $\Cat(\X, \mathbf C_{\ge
0}(\A))[\mathcal{S}^{-1}]$.
\end{corollary}

\begin{remark}
In \cite{GNPR2} we have presented some variations of the acyclic
models theorem in the monoidal and the symmetric monoidal
settings. They can also be deduced from a convenient
Cartan-Eilenberg structure.
\end{remark}

\end{large}
\end{document}